\documentclass[graybox]{svmult}

\usepackage{mathptmx}       
\usepackage{helvet}         
\usepackage{courier}        
\usepackage{type1cm}        
\usepackage{makeidx}         
\usepackage{graphicx}        
\usepackage{multicol}        
\usepackage[bottom]{footmisc}

\usepackage{amsmath, amssymb}

\makeindex             
                       
\newcommand{\field}{\mathbf{k}}
\DeclareMathOperator{\lk}{lk}
\DeclareMathOperator{\Init}{in}
\DeclareMathOperator{\st}{st}
\DeclareMathOperator{\Skel}{Skel}
\DeclareMathOperator{\Char}{char}
\DeclareMathOperator{\Soc}{Soc}
\DeclareMathOperator{\Hom}{Hom}
\DeclareMathOperator{\Span}{Span}
\DeclareMathOperator{\GL}{GL}
\DeclareMathOperator{\Gin}{Gin}

\begin{document}

\title*{Face enumeration on simplicial complexes}
\titlerunning{Face enumeration on simplicial complexes}
\author{Steven Klee and Isabella Novik}
\institute{Steven Klee \at Seattle University, Department of Mathematics, 901 12th Avenue, Seattle, WA 98122, \email{klees@seattleu.edu}
\and Isabella Novik \at University of Washington, Department of Mathematics, Box 354350, Seattle, WA 98195,  \email{novik@math.washington.edu}. \\ (Novik's research is partially supported by NSF grant DMS-1361423.)}
%
%
\maketitle

\abstract*{Let $M$ be a closed triangulable manifold, and let $\Delta$ be a triangulation of $M$. What is the smallest number of vertices that $\Delta$ can have? How big or small can the number of edges of $\Delta$ be as a function of the number of vertices? More generally, what are the possible face numbers ($f$-numbers, for short) that $\Delta$ can have? In other words, what restrictions does the topology of $M$ place on the possible $f$-numbers of triangulations of $M$? 
To make things even more interesting, we can add some combinatorial conditions on the triangulations we are considering (e.g., flagness, balancedness, etc.) and ask what additional restrictions  these combinatorial conditions impose. While only a few theorems in this area of combinatorics were known a couple of decades ago, in the last ten years or so, the field simply exploded with new results and ideas. Thus we feel that a survey paper is long overdue. As new theorems are being proved while we are typing this chapter, and as we have only a limited number of pages, we apologize in advance to our friends and colleagues, some of whose results will not get mentioned here.
}


\section{Introduction}
Let $M$ be a closed triangulable manifold, and let $\Delta$ be a triangulation of $M$. What is the smallest number of vertices that $\Delta$ can have? How big or small can the number of edges of $\Delta$ be as a function of the number of vertices? More generally, what are the possible face numbers ($f$-numbers, for short) that $\Delta$ can have? In other words, what restrictions does the topology of $M$ place on the possible $f$-numbers of triangulations of $M$? 

To make things even more interesting, we can add some combinatorial conditions on the triangulations we are considering (e.g., flagness, balancedness, etc.) and ask what additional restrictions  these combinatorial conditions impose. While only a few theorems in this area of combinatorics were known a couple of decades ago, in the last ten years or so, the field simply exploded with new results and ideas. Thus we feel that a survey paper is long overdue. As new theorems are being proved while we are typing this chapter, and as we have only a limited number of pages, we apologize in advance to our friends and colleagues, some of whose results will not get mentioned here.

The paper is structured as follows. In Section \ref{sect:definitions} we recall basic definitions pertaining to simplicial complexes. In Section \ref{Classical Results} we review classical results on the $f$-numbers of polytopes and spheres. In Section \ref{sect:manifolds} we survey many recent results on the $f$-vectors of triangulations of manifolds and pseudomanifolds. Section \ref{balanced-and-flag} is devoted to the face numbers of flag and balanced complexes. Finally, in Section \ref{sect:tools} we discuss some of the tools and techniques used in the proofs.


\section{Definitions}  \label{sect:definitions}

A \textbf{simplicial complex} $\Delta$ on a (finite) vertex set $V = V(\Delta)$ is a collection of subsets of $V(\Delta)$ that is closed under inclusion.  The elements $F \in \Delta$ are called \textbf{faces} and the maximal faces of $\Delta$ under inclusion are called \textbf{facets}.  We say that $\Delta$ is \textbf{pure} if all of its facets have the same cardinality.  The \textbf{dimension} of a face $F \in \Delta$ is $\dim(F) = |F|-1$ and the dimension of $\Delta$ is $\dim(\Delta) = \max\{\dim(F)\ : \ F \in \Delta\}$; we refer to $i$-dimensional faces as \textbf{$i$-faces}. If $\Delta$ is a $(d-1)$-dimensional simplicial complex and $k \leq d-1$, the \textbf{$k$-skeleton} of $\Delta$ is $\Skel_k(\Delta):= \{F \in \Delta \ : \ \dim(F) \leq k\}$.  

Although simplicial complexes are defined as purely combinatorial objects, each simplicial complex $\Delta$ admits a geometric realization $\|\Delta\|$ that contains a geometric $i$-simplex for each $i$-face of $\Delta$.  We typically do not distinguish between the combinatorial object $\Delta$ and the geometric object $\|\Delta\|$ and often say that $\Delta$ has certain geometric or topological properties in addition to certain combinatorial properties.

Let $\Delta$ be a $(d-1)$-dimensional simplicial complex. The main object of our study is the \textbf{$f$-vector} of $\Delta$, $f(\Delta):=(f_{-1}(\Delta), f_0(\Delta),\ldots, f_{d-1}(\Delta))$, where $f_i=f_i(\Delta)$ denotes the number of $i$-faces of $\Delta$; the numbers $f_i$ are called the \textbf{$f$-numbers} of $\Delta$. 

For several algebraic reasons, it is often more natural to study a certain (invertible) integer transformation of the $f$-vector called the \textbf{$h$-vector} of $\Delta$, $h(\Delta)=(h_0(\Delta),h_1(\Delta),\ldots,h_d(\Delta))$; it is defined by the following polynomial relation:
$$
\sum_{j=0}^d h_j(\Delta) \cdot t^{d-j}=\sum_{j=0}^d f_{j-1}(\Delta) \cdot (t-1)^{d-j}.
$$
Thus, for all $0 \leq j \leq d$,
$$h_j(\Delta)= \sum_{i=0}^j (-1)^{j-i} \binom{d-i}{d-j} f_{i-1}(\Delta) \quad \mbox{and} \quad f_{i-1}(\Delta) = \sum_{j=0}^i \binom{d-j}{d-i}h_j(\Delta).$$ 
In particular, knowing the $f$-numbers of $\Delta$ is equivalent to knowing its $h$-numbers. Moreover, since the $f$-numbers are \textit{nonnegative} integer combinations of the $h$-numbers, any upper/lower bounds on the $h$-numbers of $\Delta$ automatically imply upper/lower bounds on the $f$-numbers of $\Delta$.  

The motivating question in this field of research is the following: given a topological manifold $M$, characterize the set of $f$- or $h$-vectors of simplicial complexes $\Delta$ whose geometric realization is homeomorphic to $M$.  This question is incredibly difficult, and there are only a few manifolds for which this question is completely answered.  Instead, we relax the question in a number of ways.  

Rather than study a manifold of given topological type, we study the family of all manifold triangulations or even more relaxed families of triangulations.  Second, rather than seek complete characterizations of the sets of $f$- or $h$-vectors of the simplicial complexes in these families, we seek restrictions that can be placed on the $f$- or $h$-vectors that connect the combinatorics of the triangulation to the underlying topology of the simplicial complex.  We begin by describing some families of simplicial complexes that extend the notion of a manifold triangulation.  

If $\Delta$ is a simplicial complex and $F \in \Delta$ is a face, the local structure of $\Delta$ around $F$ is described by the (closed) \textbf{star} of $F$ in $\Delta$: $\st_{\Delta}(F):= \{G \in \Delta\ : \ F \cup G \in \Delta\}$.  Morally, when $\Delta$ triangulates a manifold, one should view the star of a face as the closed ball of faces that surround that face (however, even for manifolds, the star of a face need not be a topological ball).  Similarly, the \textbf{link} of $F$ in $\Delta$ is defined as
$$\lk_{\Delta}(F):= \{G \in \st_\Delta(F)\ : F \cap G = \emptyset\}.$$  When $F = \{v\}$ consists of a single vertex, we write $\lk_{\Delta}(v)$ in place of $\lk_{\Delta}(\{v\})$.  

If $\field$ is a field (or the ring of integers), 
the \textbf{(reduced) Betti numbers} of $\Delta$ are $\beta_i(\Delta;\field):= \dim_{\field}\widetilde{H}_i(\Delta;\field)$.  
Here $\widetilde{H}_i(-; \field)$ stands for the $i$-th simplicial reduced homology computed with coefficients in $\field$. 
The \textbf{(reduced) Euler characteristic} of $\Delta$ is 
$\widetilde{\chi}(\Delta) = \sum_{i=0}^d (-1)^{i-1}f_{i-1}(\Delta) = \sum_{i=0}^{d-1}(-1)^i\beta_i(\Delta;\field)$. Below we denote by $\mathbb{S}^j$ and $\mathbb{B}^j$ the $j$-dimensional sphere and ball, respectively.

On the one hand, we can relax the structure of a manifold triangulation (also called a \textbf{simplicial manifold}) by studying the families of simplicial homology manifolds and simplicial pseudomanifolds.  A $(d-1)$-dimensional simplicial complex $\Delta$ is a \textbf{$\field$-homology sphere} if $\widetilde{H}_*(\lk_{\Delta}(F);\field) \cong \widetilde{H}_*(\mathbb{S}^{d-|F|-1};\field)$ for all faces $F \in \Delta$ (including $F = \emptyset$). Similarly, $\Delta$ is a \textbf{$\field$-homology manifold} (without boundary) if for all nonempty faces $F\in\Delta$, the link of $F$ is a $\field$-homology $(d-|F|-1)$-sphere. 

Homology manifolds with boundary are defined in an analogous way: a simplicial complex $\Delta$ is a \textbf{$\field$-homology manifold with boundary} if (1) the link of each non-empty face $F$ of $\Delta$ has the homology of either $\mathbb{S}^{d-|F|-1}$ or $\mathbb{B}^{d-|F|-1}$, and (2) the set of all \textbf{boundary faces}, that is,
\[
\partial\Delta:= \left\{F\in\Delta \ : \ \widetilde{H}_*(\lk_{\Delta}(F);\field) \cong \widetilde{H}_*(\mathbb{B}^{d-|F|-1};\field)\right\} \cup \{\emptyset\},
\]
is a $\field$-homology $(d-2)$-manifold without boundary. 

Why do these definitions extend the definitions of manifold triangulations? Applications of the standard homology axioms (see \cite[Lemma 3.3]{Munkres-84}) show that if $F \in \Delta$ is nonempty and $p$ is a point in the relative interior of $F$, then $H_j(\|\Delta\|, \|\Delta\|-p;\field) \cong \widetilde{H}_{j-|F|}(\lk_{\Delta}(F);\field)$.  When $\|\Delta\|$ is a $(d-1)$-manifold, $H_j(\|\Delta\|, \|\Delta\|-p;\field)$ is nontrivial (and isomorphic to $\field$) if and only if $p$ is an interior point and $j = d-1$.

A $(d-1)$-dimensional simplicial complex $\Delta$ is a \textbf{normal pseudomanifold} if (1) it is pure, (2) each $(d-2)$-face (or \textbf{ridge}) of $\Delta$ is contained in exactly two facets of $\Delta$, and (3) the link of each nonempty face of dimension at most $d-3$ is connected.  In particular, every homology manifold without boundary is a normal pseudomanifold.  

On the other hand, we can impose extra structure on our simplicial complexes by studying the families of simplicial polytopes and piecewise linear (PL) manifold triangulations.  A \textbf{$d$-polytope} is the convex hull of a finite set of points in $\mathbb{R}^d$ that affinely span $\mathbb{R}^d$.  A $d$-polytope is simplicial if all of its (proper) faces are simplices.  While the boundary of any simplicial $d$-polytope is a simplicial $(d-1)$-sphere, not every simplicial sphere can be realized as the boundary of a polytope.  In fact, for $d\geq 4$, almost all simplicial $(d-1)$-spheres are not polytopal, see \cite{Kalai-88, Pfeifle-Ziegler, Nevo-Santos-Wilson}.  

One feature of a simplicial $d$-polytope is that all of its faces are realized linearly in $\mathbb{R}^d$.  A simplicial manifold $\Delta$ is a \textbf{piecewise linear} (PL) manifold if the link of each nonempty face of $\Delta$ is PL-homeomorphic to a simplex or to the boundary of a simplex of the appropriate dimension. 

The canonical example of a simplicial complex that distinguishes all of these definitions is the Poincar\'e 3-sphere.  If $\Delta$ is a triangulation of the Poincar\'e 3-sphere, then $\Delta$ is a $\mathbb{Z}$-homology 3-sphere, but not a topological 3-sphere.  Cannon \cite{Cannon-79} showed that the double suspension of any $\mathbb{Z}$-homology sphere is a topological sphere, and hence $\Sigma^2(\Delta)$, the double suspension of $\Delta$, is a triangulation of $\mathbb{S}^5$.  However, if $e$ is an edge in $\Sigma^2(\Delta)$ connecting suspension vertices, then $\lk_{\Sigma^2(\Delta)}(e) = \Delta$ is not homeomorphic to $\mathbb{S}^3$.  Therefore, $\Sigma^2(\Delta)$ is not a PL manifold.  This also shows that the link of any face in a sphere triangulation has the homology of a sphere, but need not be homeomorphic to a sphere.

Throughout this paper, $\Delta$ denotes a simplicial complex; often, it is a homology manifold. Unless we say that $\Delta$ is a homology manifold with boundary, we assume that $\partial\Delta = \emptyset$. Also, by the $f$- and $h$-numbers of a simplicial polytope $P$ we mean the  $f$- and $h$-numbers of the boundary of $P$.

\section{Classical results for polytopes and spheres} \label{Classical Results}

One of the first classical results on the structure of $f$-vectors of triangulations of manifolds is known as Klee's Dehn-Sommerville equations \cite{Klee-64}.  

\begin{theorem}\label{Dehn-Sommerville-Equations} {\rm \cite{Klee-64}\;}
Let $\Delta$ be a $(d-1)$-dimensional homology manifold (or more generally, any semi-Eulerian simplicial complex).  Then 
$$h_{d-j}(\Delta) - h_j(\Delta) = (-1)^j \binom{d}{j} \left[\widetilde{\chi}(\Delta) - \widetilde{\chi}(\mathbb{S}^{d-1})\right], \quad \mbox{for all $0 \leq j \leq d$}.$$
\end{theorem}

Thus, the $h$-vector of a homology sphere is symmetric, and hence determined by $h_0, h_1, \ldots, h_{\lfloor d/2\rfloor}$.  For a $(d-1)$-dimensional simplicial complex $\Delta$, $h_0(\Delta)=1$, while $h_d(\Delta) = \sum_{i=0}^d (-1)^{d-i}f_{i-1}(\Delta) = (-1)^{d-1} \widetilde{\chi}(\Delta)$.  Hence the Dehn-Sommerville equations provide a vast generalization of the Euler-Poincar\'e formula.

Many early results in the study of face enumeration focused on the properties of $f$-vectors of (simplicial) polytopes.  A first natural question is to ask how large or how small the $f$-numbers of a simplicial polytope could be as a function of its dimension and number of vertices.  

For the question of lower bounds, Barnette \cite{Barnette-73} showed that among all simplicial $d$-polytopes with $n$ vertices, \textbf{stacked polytopes} minimize the number of faces of each dimension.  A stacked $d$-polytope on $n$ vertices, denoted $ST(n,d)$, is defined inductively as follows.  The polytope $ST(d+1,d)$ is a $d$-simplex.  For $n>d+1$, $ST(n,d)$ is constructed by gluing a $d$-simplex to $ST(n-1,d)$ along a common facet.  

\begin{theorem}\label{LBT}  {\rm (The Lower Bound Theorem)\;} Let $\Delta$ be a $(d-1)$-dimensional connected homology manifold on $n$ vertices with $d \geq 3$.  
\begin{enumerate}
\item {\rm \cite{Barnette-73}} Then $h_2(\Delta) \geq h_1(\Delta)$.  It follows that $f_{i-1}(\Delta) \geq ST(n,d)$ for all $i$.  
\item {\rm \cite{Kalai-87}} Moreover, if $d\geq 4$, then $h_2(\Delta) = h_1(\Delta)$ if and only if $\Delta$ is the boundary complex of a stacked $d$-polytope.
\end{enumerate}
\end{theorem}
Part (2) of the Lower Bound Theorem is due to Kalai \cite{Kalai-87}; Kalai also found an alternative proof of part (1). This theorem was later generalized to all normal pseudomanifolds by Fogelsanger \cite{Fogelsanger-88} and Tay \cite{Tay-95}.  The condition that $d \geq 4$ in part (2) of the Lower Bound Theorem is necessary since $h_2 = h_1$ for all $2$-spheres, but not all $2$-spheres are stacked (e.g., the boundary of the octahedron is not stacked).

For the question of upper bounds, McMullen \cite{McMullen-70} showed that the \textbf{cyclic polytope} is a simplicial polytope with the maximal number of faces of each dimension.  A cyclic $d$-polytope on $n$ vertices, denoted $\mathcal{C}_d(n)$, is defined by taking the convex hull of $n$ distinct points on the moment curve $\varphi(t):=(t,t^2,\ldots,t^d) \subseteq \mathbb{R}^d$.  

\begin{theorem}\label{UBT} {\rm (The Upper Bound Theorem \cite{McMullen-70})\;}
Let $P$ be a simplicial $d$-polytope on $n$ vertices.  Then $h_j(P) \leq h_j(\mathcal{C}_d(n))$ for all $j$.
\end{theorem}

Stanley \cite{Stanley-75} extended the Upper Bound Theorem to all simplicial homology spheres, and Novik \cite{Novik-98} extended it to all odd-dimensional  simplicial manifolds and certain even-dimensional simplicial manifolds (see Section \ref{subs:UBT-manifolds} for more details). 

A complete characterization of $h$-vectors of simplicial polytopes was given in the celebrated $g$-theorem of Stanley \cite{Stanley-80} (necessity), and Billera and Lee \cite{Billera-Lee-81} (sufficiency). 

\begin{theorem} \label{g-theorem} {\rm (The $g$-theorem \cite{Stanley-80, Billera-Lee-81})\;} A vector $\mathbf{h} = (h_0, h_1, \ldots, h_d) \in \mathbb{Z}^{d+1}$ is the $h$-vector of a simplicial $d$-polytope if and only if 
\begin{enumerate}
\item $h_j = h_{d-j}$ for all $j$, 
\item $1 = h_0 \leq h_1 \leq \cdots \leq h_{\lfloor d/2 \rfloor}$, and
\item the numbers $(g_0, g_1,\ldots, g_{\lfloor d/2 \rfloor})$, where $g_j:=h_j-h_{j-1}$,  form an $M$-sequence.
\end{enumerate}
\end{theorem}

See Section \ref{subsection:SR-ring} for the definition of an $M$-sequence.  The $g$-theorem shows not only that the $h$-numbers of a simplicial polytope are unimodal, but also that their rate of growth is bounded. We also mention that Lee's interpretation \cite{Lee-94} of Kalai's proof of the Lower Bound Theorem implies that the portion $(g_0(\Delta),g_1(\Delta),g_2(\Delta))$ of the $g$-vector is an $M$-sequence for any connected normal pseudomanifold $\Delta$.

The unimodality of the $h$-numbers of simplicial polytopes, together with the treatment of the cases of equality, is given by the Generalized Lower Bound Theorem (GLBT, for short):

\begin{theorem}\label{GLBT} {\rm (The GLBT \cite{McMullen-Walkup-71, Stanley-80, Murai-Nevo-13})\;} 
Let $P$ be a simplicial $d$-polytope.  Then
\begin{enumerate}
\item $h_0(P) \leq h_1(P) \leq \cdots \leq h_{\lfloor d/2 \rfloor}(P)$, and
\item $h_r(P) = h_{r-1}(P)$ for some $r \leq \lfloor d/2 \rfloor $ if and only if $P$ is $(r-1)$-stacked.
\end{enumerate}

\end{theorem}
A simplicial $d$-polytope $P$ is called \textbf{$(r-1)$-stacked} if there exists a triangulation, $\Delta$, of $P$ such that $\partial\Delta = \partial P$ and $\Skel_{d-r}(\Delta) = \Skel_{d-r}(P)$.  In other words, $\Delta$ extends the triangulation of the boundary of $P$ to the interior of $P$ without introducing any interior faces of dimension  $\leq d-r$.  When $r=2$, this says that $P$ admits a triangulation with all interior faces having dimension $d-1$ or $d$.  This is precisely what it means for $P$ to be stacked in the sense of the Lower Bound Theorem (Theorem \ref{LBT}). 

\section{Face enumeration for manifolds}  \label{sect:manifolds}

In this section, we discuss extensions of the results from Section \ref{Classical Results} from spheres and polytopes to arbitrary manifolds, both with and without boundary.  Here we focus on presenting as many results as possible; the tools and techniques that are used will be presented later in the paper.

\subsection{Modified $h$-numbers and the manifold $g$-conjecture}

The Dehn-Sommerville equations (Theorem \ref{Dehn-Sommerville-Equations}) and the Generalized Lower Bound Theorem (Theorem \ref{GLBT}) show that the $h$-numbers of a simplicial polytope are symmetric, nonnegative, and unimodal.  In contrast, the $h$-numbers of a simplicial manifold need not have any of these properties -- in particular, they may be negative.  For this reason, we modify the definition of the $h$-numbers  in the following way. 

\begin{definition} \label{h-doubleprime-nums}
Let $\Delta$ be a $(d-1)$-dimensional simplicial complex and $\field$ a field.  The \textbf{$h''$-numbers} of $\Delta$ are defined as
$$
h''_j(\Delta) = 
\begin{cases}
h_j(\Delta) - \binom{d}{j} \sum_{i=0}^j (-1)^{j-i} \beta_{i-1}(\Delta;\field) & \text{ if } 0 \leq j < d, \\
\beta_{d-1}(\Delta;\field) & \text{ if } j = d.
\end{cases}
$$
\end{definition}

Note that if $\Delta$ is a $\field$-homology sphere, then $h''_j(\Delta) = h_j(\Delta)$ for all $j$.  The following result, which generalizes the Dehn-Sommerville equations and provides evidence towards a manifold $g$-theorem, further indicates that the $h''$-numbers are the ``correct" generalization of the $h$-numbers for manifolds. 

\begin{theorem} \label{h-double-prime}
Let $\Delta$ be a $(d-1)$-dimensional connected, orientable, $\field$-homology manifold.  Then
\begin{enumerate}
\item {\rm (\cite[Lemma 5.1]{Novik-98})\;} $h''_j(\Delta) = h''_{d-j}(\Delta)$ for all $0 \leq j \leq d$.
\item {\rm (\cite[Theorem 3.5]{Novik-Swartz-09:Socles})\;} $h''_j(\Delta) \geq 0$ for all $0 \leq j \leq d$.  
\item {\rm (\cite[Theorem 3.2]{Novik-Swartz-09:Gorenstein}, \cite[Theorem 2.4]{Swartz-14})\;} If all but (at most) $d+1$ of the vertex links in $\Delta$ have the weak Lefschetz property over $\field$, then 
\begin{enumerate}
\item $h''_0(\Delta) \leq h''_1(\Delta) \leq \cdots \leq h''_{\lfloor d/2 \rfloor}(\Delta)$ and
\item the numbers $(g''_j(\Delta) := h''_j(\Delta) - h''_{j-1}(\Delta))_{j=0}^{\lfloor d/2 \rfloor}$ form an $M$-sequence.
\end{enumerate}
\end{enumerate}
\end{theorem}

See Section \ref{sec:tools-lefschetz} for the basics on the weak Lefschetz property (WLP, for short). Specifically, Theorem \ref{strongLefschetz} below asserts that the boundary of any simplicial polytope has the WLP over $\mathbb{Q}$. Thus, part (3) of Theorem \ref{h-double-prime}  applies  if $\field=\mathbb{Q}$ and most of vertex links of $\Delta$ are polytopal.

If $\Delta$ is disconnected with connected components $\Delta_1, \ldots, \Delta_s$, it is easy to show that $h''_0(\Delta) = 1$ and $h''_j(\Delta) = \sum_{i=1}^s h''_j(\Delta_i) $ for $1 \leq j \leq d$.  Therefore parts (2) and (3) of Theorem \ref{h-double-prime} continue to hold for disconnected manifolds, and part (1) holds as long as $0 < j < d$. Further, part (2) of Theorem \ref{h-double-prime} continues to hold for non-orientable manifolds and manifolds with boundary.  

The celebrated \textbf{$g$-conjecture for spheres} posits that the conditions of the $g$-theorem (Theorem \ref{g-theorem}) hold for {\em all} simplicial (homology) spheres.  Kalai proposed a far-reaching generalization of this conjecture: he conjectured that part (3) of Theorem~\ref{h-double-prime} continues to hold for {\em all} orientable, $\field$-homology manifolds.  We refer to this latter conjecture as the \textbf{manifold $g$-conjecture}. 

Part (3) of Theorem~\ref{h-double-prime} shows that if $\field$-homology spheres have the WLP, then the manifold $g$-conjecture also holds. In general, very little is known about either version of the $g$-conjecture.  Kubitzke and Nevo \cite{Kubitzke-Nevo-09}, building on work of Brenti and Welker \cite{Brenti-Welker-08}, proved that the barycentric subdivision of a simplicial homology sphere satisfies the $g$-conjecture for spheres.  Murai \cite{Murai-10:Barycentric} extended this result to show that the barycentric subdivision of a simplicial homology manifold satisfies the manifold $g$-conjecture.  Some other classes of spheres for which the $g$-conjecture is known to hold include spheres with few vertices or few edges \cite[Theorem 4.1]{Swartz-14}, as well as PL-spheres that can be obtained from the boundary of a simplex by certain (in fact, most) bistellar flips \cite[Theorem 3.1, Theorem 3.2]{Swartz-14}.

\subsection{The Upper Bound Theorem for Manifolds}  \label{subs:UBT-manifolds}

In 1957, Motzkin proposed the Upper Bound Conjecture (UBC), which states that a cyclic $d$-polytope on $n$ vertices has the componentwise maximal $f$-vector among all $d$-polytopes on $n$ vertices.  We use $\mathcal{C}_d(n)$ to denote a cyclic $d$-polytope on $n$ vertices.  Since a small perturbation of the vertices of a polytope will (1) result in a simplicial polytope and (2) only increase the number of faces,  it suffices to verify Motzkin's conjecture for simplicial polytopes.  

In 1964, V.~Klee \cite{Klee-64:vertices} extended the UBC, proposing that it continues to hold for all \textbf{Eulerian} simplicial complexes; i.e, those simplicial complexes $\Delta$ with the property that $\widetilde{\chi}(\lk_{\Delta}(F)) = (-1)^{d-|F|-1}$ for all $F \in \Delta$ (including the empty face).  Moreover, he proved the UBC for all Eulerian complexes with sufficiently many vertices.

McMullen \cite{McMullen-70} proved the UBC for polytopes in 1970, and Stanley \cite{Stanley-75} proved the UBC for all simplicial homology spheres in 1975\footnote{We further refer to Stanley's recent article \cite{Stanley-14} which gives his personal account of how he came to the proof of the Upper Bound Theorem.}.  Novik \cite{Novik-98} extended these results further to homology manifolds that are also Eulerian simplicial complexes.

\begin{theorem}  \label{UBT-extensions}
Let $\Delta$ be an orientable $\field$-homology manifold of dimension $d-1$ with $n$ vertices.  If either
\begin{enumerate}
\item {\rm (\cite[Theorem 1.4]{Novik-98})} $\; d-1$ is odd, or
\item {\rm (\cite[Theorem 1.5]{Novik-98}, \cite[Theorem 5.3]{Novik-05})}  $\; d-1=2k$ is even and 
\begin{equation} \label{UBC-homology-ineq} 
\beta_k(\Delta;\field) \leq 2\beta_{k-1} + 2\sum_{i=0}^{k-3} \beta_i(\Delta;\field),
\end{equation}
\end{enumerate}
then $h_j(\Delta) \leq h_j(\mathcal{C}_d(n))$ for all $j$.  In particular, $f_{i-1}(\Delta) \leq f_{i-1}(\mathcal{C}_d(n))$ for all $i$.
\end{theorem}

Therefore, the Upper Bound Theorem (Theorem \ref{UBT}) continues to hold for all simplicial odd-dimensional manifolds and for any even-dimensional manifolds for which \eqref{UBC-homology-ineq} holds over $\mathbb{Z}/2\mathbb{Z}$.  When $\dim(\Delta) = 2k$ is even, \eqref{UBC-homology-ineq} is satisfied by all $\field$-homology manifolds for which either $\beta_k(\Delta; \field) = 0$ or $(-1)^k(\widetilde{\chi}(\Delta)-1) \leq 0$. 

It is worth adding that the ``in particular" part (i.e., the statement on $f$-vectors) of Theorem \ref{UBT-extensions} was extended to all odd-dimensional Eulerian pseudomanifolds with isolated singularities in \cite{Hersh-Novik-02}; moreover, Novik and Swartz \cite{Novik-Swartz-12}, proved the UBT for all Eulerian pseudomanifolds with so-called homologically isolated singularities as long as $n\geq 3d-4$. However, Klee's conjecture \cite{Klee-64:vertices} that the UBC holds for all Eulerian complexes remains wide open.

\subsection{K\"uhnel's conjectures}
By  Poincar\'e duality, the Euler characteristic of any odd-dimensional simplicial manifold is equal to zero. Motivated by the UBT, K\"uhnel \cite{Kuhnel-90} conjectured an upper bound on the Euler characteristic of even-dimensional PL manifolds.  He also proved this conjecture  for all such manifolds with sufficiently many vertices (see \cite[p.~66]{Kuhnel-95}). Novik \cite{Novik-98} verified K\"uhnel's conjecture for all $2k$-dimensional simplicial manifolds with $n\geq 4k+3$ or $n\leq 3k+3$ vertices, and Novik and Swartz \cite{Novik-Swartz-09:Socles} proved the conjecture in full generality:

\begin{theorem} \label{Kuhnel-1} {\rm \cite[Theorem 4.4]{Novik-Swartz-09:Socles}\;} Let $\Delta$ be a $2k$-dimensional orientable $\field$-homology manifold with $n$ vertices. Then
$$
(-1)^k \binom{2k+1}{k}\left(\widetilde{\chi}(\Delta)-1\right)\leq \binom{n-k-2}{k+1}.
$$
Moreover, equality holds if and only if $\Delta$ is $(k+1)$-neighborly.
\end{theorem}
A simplicial complex $\Delta$ is \textbf{$s$-neighborly} if every $s$-vertex subset of $V(\Delta)$ is a face of $\Delta$. For instance, the cyclic polytope $\mathcal{C}_d(n)$ is $\lfloor d/2\rfloor$-neighborly. 

K\"uhnel also conjectured upper bounds on the individual Betti numbers; these conjectures are recorded in \cite[Conjecture 18]{Lutz-05}. To this end, Murai \cite{Murai-15} showed:

\begin{theorem} \label{Kuhnel-2} Let $\Delta$ be a $(d-1)$-dimensional $\field$-homology manifold with $n$ vertices.
\begin{enumerate}
\item {\rm \cite[Theorem 5.1]{Murai-15} \;} If $d=2k+1$, then $\binom{2k+1}{k}\beta_k(\Delta;\field)\leq \binom{n-k-2}{k+1}$.
\item {\rm \cite[Theorem 5.2]{Murai-15} \;} If all vertex links of $\Delta$ are polytopal, then 
\begin{equation} \label{eq:bound-on-beta}
\binom{d+1}{j+1}\beta_j(\Delta;\mathbb{Q})\leq \binom{n-d+j-1}{j+1} \quad \mbox{for all } 0\leq j \leq \lfloor d/2\rfloor -1.
\end{equation}
\end{enumerate}
\end{theorem}

The same result, but restricted to the class of orientable homology manifolds, was proved in \cite{Novik-Swartz-09:DS}. 

Note that if $\beta_j(\Delta;\mathbb{Q})\neq 0$, then equation \eqref{eq:bound-on-beta} implies that $$\binom{d+1}{j+1}\leq \binom{n-d+j-1}{j+1}, \mbox{ hence } d+1\leq n-d+j-1, \mbox{ and so } n\geq 2d+2-j.$$ In fact, for PL-manifolds, Brehm and K\"uhnel \cite{Brehm-Kuhnel-87} proved the following stronger result: if $n=2d+1-j$, then the homotopy groups $\pi_0(\Delta), \pi_1(\Delta), \ldots, \pi_j(\Delta)$ are all trivial. Their proof relies on PL Morse theory.

\subsection{The Lower Bound Theorem for manifolds} 
According to the Lower Bound Theorem (Theorem \ref{LBT}), any connected homology manifold (or even a normal pseudomanifold) $\Delta$ of dimension $d-1\geq 2$ satisfies $g_2(\Delta)\geq 0$. In an effort to understand how the topology of $\|\Delta\|$ affects the combinatorics of $\Delta$, Kalai \cite{Kalai-87} conjectured that for manifolds of dimension $d-1\geq 3$ this inequality can be strengthened to $g_2(\Delta)\geq \binom{d+1}{2}\beta_1(\Delta;\field)$.  

\begin{theorem} \label{LBT-Murai} {\rm \cite[Theorem 5.3]{Murai-15}\;} Let $\Delta$ be a normal pseudomanifold of dimension $d-1\geq 3$. Then 
\begin{enumerate}
\item $g_2(\Delta)\geq \binom{d+1}{2}\left(\beta_1(\Delta;\field)-\beta_0(\Delta;\field)\right)$.
\item Moreover, equality holds if and only if $\Delta$ is a stacked simplicial manifold.
\end{enumerate}
\end{theorem}

Note that Theorem \ref{LBT-Murai} strengthens part (2) of Theorem \ref{Kuhnel-2} in the case that $j=1$. Theorem \ref{LBT-Murai} is also due to Murai \cite{Murai-15}. In the special case that $\Delta$ is an orientable  $\field$-homology manifold, part (1) of this result, along with the $d\geq 5$ case of part (2), was proved in \cite[Theorem 5.2]{Novik-Swartz-09:Socles}. The case of equality for 3-dimensional homology manifolds (both orientable and non-orientable) is due to Bagchi \cite{Bagchi:mu-vector}. 

We will discuss \textbf{stacked manifolds}, and more generally, \textbf{$r$-stacked manifolds} in Section \ref{section:r-stacked}. For now, it suffices to say that stacked manifolds form a natural generalization of \textbf{stacked spheres} -- the boundary complexes of stacked polytopes. In particular, for $d\geq 5$, a $(d-1)$-dimensional simplicial manifold is stacked if and only if all its vertex links are stacked spheres.

\subsection{The $\sigma$- and $\mu$-vectors}  \label{section:sigma-and-mu}

Bagchi and Datta \cite{Bagchi-Datta-13, Bagchi:mu-vector}, motivated by  Morse theory (see \cite[Remark 2.11]{Bagchi-Datta-13}), introduced the notion of $\mu$-numbers. The following definition is Murai's \cite{Murai-15} slight modification of their original definition. Denote by $\Delta_W$ the subcomplex of $\Delta$ induced by vertices in $W\subseteq V(\Delta)$ and recall that $\beta_0(\{\emptyset\};\field)=0$, while $\beta_{-1}(\{\emptyset\};\field)=1$.

\begin{definition} Let $\Delta$ be a $(d-1)$-dimensional simplicial complex and $\field$ a field. The $\sigma$- and $\mu$-vectors of $\Delta$, $(\sigma_{-1}, \sigma_0,\ldots, \sigma_{d-1})$ and $(\mu_0,\ldots,\mu_{d-1})$, are defined by
\[
\sigma_j=\sigma_j(\Delta;\field):=\sum_{W\subseteq V} \frac{\beta_j(\Delta_W;\field)}{\binom{|V|}{|W|}} \mbox{  and }
\mu_j=\mu_{j}(\Delta;\field):=\sum_{v\in V} \frac{\sigma_{j-1}(\lk_\Delta(v);\field)}{f_0(\lk_\Delta(v))+1}.
\]
\end{definition}

A few remarks are in order. First, note that the $\mu$-numbers are rational, but not necessarily integer numbers. An important property of the $\mu$-numbers is that they satisfy the following inequalities, just as the usual Morse numbers do: 
\begin{theorem} {\rm \cite{Bagchi:mu-vector}\;} \label{mu-vs-beta} Let $\Delta$ be any simplicial complex. Then for all $j$,
\begin{enumerate}
\item $\mu_j(\Delta;\field)\geq \beta_j(\Delta;\field)$ if $j>0$, and $\mu_0 (\Delta;\field)\geq \beta_0(\Delta;\field)+1$.
\item $\sum_{k=0}^j (-1)^{j-k} \mu_k(\Delta;\field)\geq \sum_{k=0}^j(-1)^{j-k} \beta_j(\Delta;\field) + (-1)^j$.
\end{enumerate}
\end{theorem}

Murai \cite{Murai-15} used Hochster's formula to interpret the $\mu$-numbers of $\Delta$ in terms of the graded Betti numbers of the Stanley-Reisner ring of $\Delta$ (see Section \ref{section:gradedBetti} for more details on the graded Betti numbers). He then used this interpretation, along with some recent results from commutative algebra, to prove Theorems \ref{Kuhnel-2} and \ref{LBT-Murai}. In fact, his proofs show that these two theorems continue to hold even if one replaces the Betti numbers of $\Delta$ with the $\mu$-numbers of $\Delta$. Specifically, one obtains the following strengthening of Theorem \ref{LBT-Murai}.

\begin{theorem} \label{LBT-Murai-mu}  Let $\Delta$ be a normal pseudomanifold of dimension $d-1\geq 3$. Then $g_2(\Delta)\geq \binom{d+1}{2}\left(\mu_1(\Delta;\field)-\mu_0(\Delta;\field)+1\right)$.
\end{theorem}

One of Kalai's long-standing conjectures posits that for any connected homology manifold $\Delta$ of dimension $d-1$, $g_2(\Delta)/\binom{d+1}{2}$ is at least as large as the minimum number of generators of the fundamental group of $\Delta$.  Currently, the only progress on this conjecture is a result of Klee \cite{Klee-09} stating that if $\Delta$ is balanced (see Section \ref{section:balanced} for a definition), then $h_2(\Delta) \geq \binom{d}{2} m(\Delta)$, where $m(\Delta)$ is the minimum number of generators of $\pi_1(\Delta)$.  Kalai's conjecture along with Theorem \ref{LBT-Murai-mu} prompted Swartz \cite{Swartz-Oberwolfach} to ask if $\mu_1(\Delta)$ provides a lower bound for the number of generators of $\pi_1(\Delta)$.

\subsection{Stronger manifold $g$-conjectures}

In view of Theorem \ref{LBT-Murai}, it is natural to define the following invariants for $r\leq\lfloor d/2\rfloor$:
\begin{eqnarray*}
\widetilde{g}_r(\Delta)&:=&g_r(\Delta)-\binom{d+1}{r}\sum_{j=1}^r (-1)^{r-j}\beta_{j-1}(\Delta;\field); \quad \mbox{equivalently,} \\
 \widetilde{g}_r(\Delta)&=&h''_r(\Delta)-h''_{r-1}(\Delta)-\binom{d}{r-1}\beta_{r-1}(\Delta;\field).
\end{eqnarray*}
In particular, $\widetilde{g}_r(\Delta)\leq g''_r(\Delta)$. The following strengthening of part (3) of Theorem \ref{h-double-prime} holds. See Section \ref{sec:tools-lefschetz} for the basics on the weak Lefschetz property.

\begin{theorem} \label{tilde-g} Let $\Delta$ be an orientable $(d-1)$-dimensional, $\field$-homology manifold and assume that all vertex links of $\Delta$ have the WLP over $\field$. Then
\begin{enumerate} 
\item The numbers $\widetilde{g}_r(\Delta)$ are nonnegative for all $r\leq d/2$.
\item Moreover, $(\widetilde{g}_0(\Delta), \widetilde{g}_1(\Delta),\ldots,\widetilde{g}_{\lfloor d/2\rfloor}(\Delta))$ is an $M$-sequence.
\end{enumerate}
\end{theorem}

A straightforward computation shows that if $\Delta$ is disconnected, then for $r\neq 0$, $\widetilde{g}_r(\Delta)$ equals the sum of the $\widetilde{g}_r$-numbers of the connected components of $\Delta$, while $\widetilde{g}_0(\Delta)=1$. Hence the above statement for disconnected manifolds follows easily from the connected case. In the connected case, part (1) of Theorem \ref{tilde-g} is implicit in \cite[Eq.~(9)]{Novik-Swartz-09:Socles}, and part (2) is due to Murai and Nevo \cite[Theorem 5.4]{Murai-Nevo-14}.

Bagchi \cite[Conjecture 2]{Bagchi:mu-vector} proposed the following strengthening of part (1) of Theorem \ref{tilde-g} and proved it for the case of $d=4$:

\begin{conjecture} \label{Bagchi-conjecture} Let $\Delta$ be a $(d-1)$-dimensional $\field$-homology manifold. Then 
\[
g_r(\Delta)\geq \binom{d+1}{r}\left[(-1)^r+\sum_{j=1}^r (-1)^{r-j}\mu_{j-1}(\Delta;\field)\right]
\quad \mbox{for all } r\leq d/2.
\]
\end{conjecture}

Note that Theorem \ref{LBT-Murai-mu} proves the $r=2$ case of this conjecture even for normal pseudomanifolds. Presently, this conjecture is wide open for $r>2$, as is a weaker version of the conjecture positing that if $\Delta$ is any $\field$-homology manifold, then  $\widetilde{g}_r(\Delta) \geq 0$ for all $r\leq d/2$.

\subsection{$r$-stacked manifolds and the cases of equality} \label{section:r-stacked}
In the previous sections we discussed many inequalities on the $f$-numbers of manifolds. The next natural question is whether these inequalities are sharp, and if so, when equality is achieved. To provide the answers we need to let manifolds with boundary enter the scene.

\begin{definition} Let $0\leq r\leq d-1$ be an integer.
\begin{enumerate} 
\item A $d$-dimensional $\field$-homology manifold with boundary $\Delta$ is \textbf{$r$-stacked} if it has no interior faces of dimension $\leq d-r-1$, that is, $\Skel_{d-r-1}(\Delta)=\Skel_{d-r-1}(\partial\Delta)$.
\item A $(d-1)$-dimensional $\field$-homology manifold without boundary is \textbf{$r$-stacked} if it is the boundary of an $r$-stacked $\field$-homology $d$-manifold with boundary.
\item A $\field$-homology manifold without boundary is \textbf{locally $r$-stacked} if all of its vertex links are $r$-stacked.
\end{enumerate}
\end{definition}

The $0$-stacked manifolds are precisely the (boundaries of) simplices; $1$-stacked manifolds are also known as stacked manifolds. The family of $r$-stacked manifolds was introduced in \cite{Murai-Nevo-14} as a natural generalization of the family of $r$-stacked polytopes (mentioned at the end of Section \ref{Classical Results}).  It is easy to see that being $r$-stacked implies being locally $r$-stacked. In fact, for small $r$, these two notions are equivalent:

\begin{theorem} \label{local-stackedness} {\rm \cite[Theorem 4.6]{Murai-Nevo-14}\;} 
For $1\leq r<d/2$, a $\field$-homology $(d-1)$-manifold without boundary is $(r-1)$-stacked if and only if it is locally $(r-1)$-stacked.
\end{theorem}

Below we summarize several results on the cases of equality. The first two are due to Murai and Nevo \cite{Murai-Nevo-14}, and the last one is due to Bagchi \cite{Bagchi:mu-vector}.

\begin{theorem} \label{h''=0} {\rm \cite[Theorem 3.1]{Murai-Nevo-14}\;}
Let $1\leq r\leq d$, and let $\Delta$ be a $\field$-homology $(d-1)$-manifold with boundary. Then $h''_r(\Delta)=0$ if and only if $\Delta$ is $(r-1)$-stacked.
\end{theorem}

\begin{theorem} \label{tilde g=0} {\rm \cite[Corollary 5.8]{Murai-Nevo-14}\;}
Let $1\leq r< d/2$, and let $\Delta$ be a $\field$-homology $(d-1)$-manifold without boundary. If all vertex links of $\Delta$ have the WLP, then $\widetilde{g}_r(\Delta)=0$ if and only if $\Delta$ is $(r-1)$-stacked.
\end{theorem}

If $d=2k$,  does $\widetilde{g}_k(\Delta)=0$ imply that $\Delta$ is $(k-1)$-stacked?

\begin{theorem} \label{g_r-equal-mu} {\rm \cite{Bagchi:mu-vector}\;} Let $1\leq r\leq d/2$ and let $\Delta$ be a locally $(r-1)$-stacked $\field$-homology $(d-1)$-manifold without boundary.  Then 
\begin{equation}  \label{Bagchi-equality} 
g_r(\Delta)= \binom{d+1}{r}\left[(-1)^r+\sum_{j=1}^r (-1)^{r-j}\mu_{j-1}(\Delta;\field)\right].
\end{equation}
\end{theorem}

Bagchi \cite{Bagchi:mu-vector} also conjectured that the converse holds, namely, if Equation \eqref{Bagchi-equality} is satisfied as equality, then $\Delta$ is locally $(r-1)$-stacked; he verified this conjecture for $3$-dimensional homology manifolds.

Modulo the results mentioned in this section, the next natural question is whether $r$-stacked manifolds exist. The $1$-stacked manifolds without boundary (of dimension $\geq 2$) are precisely the elements of the Walkup's class introduced in \cite{Walkup-70} (see also \cite{Datta-Murai}). Each such manifold is obtained by starting with several disjoint boundary complexes of the $d$-simplex and repeatedly forming connected sums and/or handle additions. Similarly, each $1$-stacked manifold with boundary (of dimension $\geq 2$) is obtained by starting with a number of $d$-simplices and repeatedly forming connected unions and/or handle additions; see \cite{Datta-Murai} for more details. Therefore, all $1$-stacked manifolds without boundary are spheres or connected sums of sphere bundles over $\mathbb{S}^1$. An analogous statement holds for $1$-stacked manifolds with boundary.

A construction of $r$-stacked manifolds with or without boundary for all $r$ and $d$ is due to 
Klee and Novik \cite{Klee-Novik-12}. Specifically, for each pair $0\leq r\leq d-1$, they construct a certain $d$-manifold with boundary, $B(r,d+1)$, that is $r$-stacked and whose boundary triangulates $\mathbb{S}^r\times \mathbb{S}^{d-1-r}$. The construction of $B(r,d+1)$ is so simple that we sketch it here. The vertex set is $X\sqcup Y$, where $X=\{x_1,\ldots, x_{d+1}\}$ and $Y=\{y_1,\ldots, y_{d+1}\}$, and the facets of $B(r,d+1)$ are all $(d+1)$-subsets of $X\sqcup Y$ of the form $\{z_1,\ldots,z_{d+1}\}$, where $z_i\in\{x_i,y_i\}$ and the number of indices $k$ such that $z_k$ and $z_{k+1}$ were drawn from different sides does not exceed $r$.

As we saw above, being $1$-stacked imposes severe restrictions on the topology of the underlying manifold. More generally, Swartz \cite{Swartz-14} proved the following:

\begin{theorem} {\rm \cite[Section 6]{Swartz-14}\;} If $\Delta$ is an $r$-stacked PL manifold with boundary, then $\Delta$ has a handle decomposition consisting of handles of index $\leq r$.
\end{theorem}
Swartz \cite{Swartz-14} also asked if every PL $d$-manifold with boundary which has a handle decomposition using only handles of index $\leq r$ admits an $r$-stacked triangulation.  While the answer is yes for $r=0,1,d-1$, for all other values of $r$ this question is wide open. We refer to \cite[Section 1.6]{Benedetti} for the definition and discussion of some properties of handle decompositions.

\subsection{The Dehn-Sommerville equations for manifolds with boundary}
Are there analogs of Dehn-Sommerville equations for manifolds with boundary? To this end, Gr\"abe \cite{Grabe-87} (see also \cite[Theorem 3.1]{Novik-Swartz-09:DS}) established the following result.

\begin{theorem} {\rm \cite{Grabe-87}\;} \label{DS-boundary}
Let $\Delta$ be a $\field$-homology $(d-1)$-manifold with boundary. Then
\[
h_{d-i}(\Delta)-h_i(\Delta)=\binom{d}{i} (-1)^{d-i-1}\widetilde{\chi}(\Delta)-g_i(\partial\Delta) \quad \mbox{for all } 0\leq i\leq d.
\]
In particular, if $\Delta$ is a $\field$-homology ball, then $h_{i}(\Delta)-h_{d-i}(\Delta)=g_i(\partial\Delta)$.
\end{theorem}

In fact, versions of the same result but in the $f$-vector form go back at least to Macdonald \cite{Macdonald-71}. The above result, along with known inequalities on the $h$-numbers of manifolds with boundary (such as the fact that $h''_i\geq 0$ for all $i$), may lead to a ``correct" version of the $h''$-numbers for manifolds with boundary. One such definition was proposed in \cite[Definition 3.3]{Novik-Swartz-09:DS}.

\subsection{The Lower Bound Theorem for manifolds with boundary}
What inequalities do the $f$-vectors of manifolds with boundary satisfy in addition to the non-negativity of $h''$-numbers? Not much is known at present. One existing result is an analog of the Lower Bound Theorem. In the following, we use $f_0^\circ(\Delta):=f_0(\Delta)-f_0(\partial\Delta)$ to denote the number of interior vertices of $\Delta$.

\begin{theorem} \label{LBT-boundary} Let $\Delta$ be a connected $(d-1)$-dimensional $\field$-homology manifold with boundary. 
\begin{enumerate} \item If $d\geq 5$, then $h_2(\Delta)\geq f_0^\circ(\Delta)+\binom{d}{2}\beta_1(\partial\Delta; \field)+d\beta_0(\partial\Delta;\field)$.
\item If $d=4$ and $\Char\field=2$, then $h_2(\Delta)\geq f_0^\circ(\Delta)+3\beta_1(\partial\Delta;\field)+4\beta_0(\partial\Delta;\field)$.
\end{enumerate}
\end{theorem}

Theorem \ref{LBT-boundary} provides a strengthening of Kalai's result \cite[Theorem 1.3]{Kalai-87} asserting that connected manifolds with boundary satisfy $h_2(\Delta)\geq f_0^\circ(\Delta)$. Theorem \ref{LBT-boundary} with an additional assumption that $\Delta$ has an orientable boundary was proved in \cite[Theorem 5.1]{Novik-Swartz-09:DS}. However, the orientability assumption was only needed when $d\geq 5$ to invoke the inequality $g_2(\partial\Delta)\geq \binom{d}{2}(\beta_1(\partial\Delta; \field)-\beta_0(\partial\Delta;\field))$ of Theorem \ref{LBT-Murai}. As Theorem~\ref{LBT-Murai} has since been established for non-orientable manifolds of dimension $\geq 3$, this additional assumption is no longer required. All 2-dimensional homology manifolds without boundary are disjoint unions of closed surfaces; hence using the field of characteristic two in part (2) of the theorem maximizes the relevant Betti numbers.

Theorem \ref{LBT-boundary} is optimal (see discussion at the end of Section 5 in \cite{Novik-Swartz-09:DS}). Thus the next question to address is the case of equality.  A conjecture to this end was proposed in \cite[Conjecture 5.7]{Novik-Swartz-09:DS}. However, there are counterexamples to this conjecture that will be discussed in some future paper, and the question of when equality holds remains wide open.

Are there any lower bounds on the higher $h$-numbers of manifolds with boundary that are analogous to those in Theorem  \ref{LBT-boundary}? At the moment we do not have even conjectural answers to this question.

\subsection{Characterizing $f$-vectors of manifolds}
At present we have only a {\em conjectural} characterization of the $f$-vectors of simplicial spheres of dimension five and above.  As such, it would be overly optimistic to hope for a complete characterization of the $f$-vectors of other manifolds of higher dimensions. However, in dimensions up to four, complete characterizations are available for a few manifolds.  These manifolds are the topic of this subsection.

For 2-manifolds, it follows from the Dehn-Sommerville equations that the entire $f$-vector is determined by the number of vertices and the Euler characteristic of a manifold in question. The $f$-vectors of \textbf{all} (closed) 2-manifolds were characterized by Ringel \cite{Ringel} for non-orientable manifolds and by Jungerman and Ringel \cite{Jungerman-Ringel-80} for orientable manifolds.

In dimension three, Walkup \cite{Walkup-70} characterized all possible $f$-vectors of $\mathbb{S}^3$, $\mathbb{R}P^3$, $\mathbb{S}^2\times\mathbb{S}^1$, and $\mathbb{S}^2\tilde{\times}\mathbb{S}^1$, the non-orientable  $\mathbb{S}^2$-bundle over $\mathbb{S}^1$ . Moreover, Lutz, Sulanke, and Swartz \cite{Lutz-Sulanke-Swartz-09} determined possible $f$-vectors of 20 additional $3$-manifolds, namely, $(\mathbb{S}^2\times\mathbb{S}^1)^{\#k}$ and $(\mathbb{S}^2\tilde{\times}\mathbb{S}^1)^{\#k}$ for $k\in\{2,3,\ldots,8,10,11,14\}$. 

In dimension four, Walkup \cite{Walkup-70} characterized all possible $f$-vectors of $\mathbb{S}^4$, while Swartz \cite{Swartz-09} characterized all possible $f$-vectors of $\mathbb{C}P^2$, $(\mathbb{S}^2\times\mathbb{S}^2)^{\#2}$, $\mathbb{S}^3\times\mathbb{S}^1$, and $K3$-surfaces. The case of the non-orientable $\mathbb{S}^3$-bundle over $\mathbb{S}^1$ was added in \cite{Chestnut-Sapir-Swartz-08}.

By the Dehn-Sommerville equations, if $M$ is a manifold of dimension $3$ or $4$ and $\Delta$ is a triangulation of $M$, then $g_1(\Delta)$ and $g_2(\Delta)$ determine the entire $f$-vector of $\Delta$. Since, in addition, $(1,g_1,g_2)$ must form an $M$-sequence, it is perhaps not so surprising that all known characterizations of the $f$-vectors for $3$- and $4$-dimensional manifolds are of the form $a \leq g_1$, $b\leq g_2\leq\binom{g_1+1}{2}$. 

What can be said about characterizing possible $f$-vectors of triangulations of manifolds with boundary? Here the situation seems even more hopeless: at present there is not even a conjecture for characterizing the set of $f$-vectors of balls of dimension six and above, see \cite{Kolins-11}. 

It follows from the $g$-conjecture for spheres and Theorem \ref{DS-boundary} that if $\Delta$ is a simplicial ball with $h$-vector $(h_0,h_1,\ldots, h_{d-1}, h_d)$, then $(h_0-h_d,h_1-h_{d-1},\ldots, h_{\lfloor d/2\rfloor}-h_{\lceil d/2\rceil})$ is an $M$-sequence, and moreover, so is
\begin{equation}  \label{Billera-Lee-ineq}
(h_0-h_{d+k}, h_1-h_{d+k-1},\ldots, h_m-h_{d+k-m}) \;\; \forall  k=0,1,\ldots, d+1, \;  m=\lfloor (d+k-1)/2\rfloor,
\end{equation}
where we take $h_i= 0$ if $i> d$. Billera and Lee \cite{Billera-Lee-81} conjectured that these conditions are also sufficient. In other words, they conjectured that every non-negative integer vector $(h_0,\ldots, h_{d-1},h_d)$ with the property that all difference vectors as in Equation \eqref{Billera-Lee-ineq} are $M$-sequences provides the $h$-vector of a simplicial ball. While Lee and Schmidt \cite{Lee-Schmidt-12} confirmed this conjecture in dimensions three and four, Kolins \cite{Kolins-11} disproved it in dimensions five and higher.

\section{Manifolds with additional structure} \label{balanced-and-flag}

In this section we discuss two combinatorial conditions that can be imposed on a simplicial complex and explore the ways in which these conditions affect the $f$- or $h$-vectors of the resulting complexes.  While these two families of simplicial complexes seem very different on the surface, the sets of $f$- and $h$-vectors arising from them are surprisingly similar and oftentimes directly related.

\subsection{Balanced complexes} \label{section:balanced}

\begin{definition}
A $(d-1)$-dimensional simplicial complex $\Delta$ is \textbf{balanced}  if its $1$-skeleton (viewed as a graph) admits a proper $d$-coloring.  
\end{definition}

Since the graph of a $(d-1)$-simplex is the complete graph on $d$ vertices, balanced complexes are the $(d-1)$-dimensional complexes with the smallest possible chromatic number.

The family of balanced simplicial complexes was introduced by Stanley \cite{Stanley-79}.  In that paper what we will call balanced complexes are called completely balanced simplicial complexes, as he allows more general colorings in which every facet could have, for example, two blue vertices, one green vertex, and three red vertices.  

The canonical example of a balanced simplicial complex is the barycentric subdivision of a simplicial complex (or a regular CW complex), $\Delta$, which is denoted as $\text{Sd}(\Delta)$.  The vertices of $\text{Sd}(\Delta)$ correspond to the faces of $\Delta$, and the faces of $\text{Sd}(\Delta)$ correspond to chains of nonempty faces in $\Delta$ under inclusion.  The barycentric subdivision of $\Delta$ is balanced by assigning color $k$ to a vertex in $\text{Sd}(\Delta)$ if its corresponding face in $\Delta$ has cardinality $k$.

\subsubsection{Upper bounds for balanced complexes}  

In general the $h$-vector of a sphere or ball is an $M$-sequence, meaning there is a family of monomials, closed under divisibility, with $h_j$ monomials of degree $j$ for all $j$ (see Section \ref{subsection:SR-ring} for further details).  The following result imposes stronger restrictions on the $h''$-vectors of balanced manifolds.  

\begin{theorem} \label{h-double-prime-balanced}
Let $\Delta$ be a balanced manifold, with or without boundary.  Then there exists a simplicial complex $\Gamma$ such that $h''(\Delta) = f(\Gamma)$.  
\end{theorem}

 In particular, it follows that for a balanced manifold, there is a family of \textit{squarefree} monomials, closed under divisibility, with $h''_j$ monomials of degree $j$ for all~$j$.  

Recall that if $\Delta$ is a simplicial (homology) sphere or ball, then $h''(\Delta) = h(\Delta)$.  In that case, Theorem \ref{h-double-prime-balanced} is due to Stanley \cite[Corollary 4.5]{Stanley-79}. For manifolds, Babson and Novik \cite[Theorem 6.6]{Babson-Novik} proved a weaker statement than the one in Theorem~\ref{h-double-prime-balanced}. However, the ideas of their proof, combined with the proof of \cite[Theorem 4.6]{Novik-Swartz-09:Socles}, imply Theorem \ref{h-double-prime-balanced}.

In view of the above result, it is natural to extend Kalai's manifold $g$-conjecture for balanced manifolds in the following way. 

\begin{conjecture}
Let $\Delta$ be a balanced, orientable $\field$-homology manifold without boundary.  Then $(g''_j(\Delta):=h''_j(\Delta)-h''_{j-1}(\Delta))_{j=0}^{\lfloor d/2 \rfloor}$ is the $f$-vector of a simplicial complex.
\end{conjecture}

Murai \cite[Theorem 1.3]{Murai-10:Barycentric} proved this conjecture in the specific case that $\Delta$ arises as the barycentric subdivision of a homology manifold with or without boundary.

It is also natural to ask if there is a balanced analog of the cyclic polytope; i.e., a balanced simplicial sphere with the componentwise maximal $f$- or $h$-numbers among all balanced spheres with a given dimension and number of vertices.

\subsubsection{Lower bounds for balanced bomplexes}

The following result extends the Lower Bound Theorem to balanced simplicial complexes.

\begin{theorem} \label{balancedLBT}
Let $\Delta$ be a balanced, connected normal pseudomanifold of dimension $d-1$ with $d \geq 3$.  Then 
\begin{equation} \label{balancedLBT-inequality}
2h_2(\Delta) \geq (d-1) h_1(\Delta).
\end{equation}
Moreover, when $d \geq 4$, equality holds in Equation \eqref{balancedLBT-inequality} if and only if $\Delta$ is a stacked cross-polytopal sphere.
\end{theorem}

Recall that a stacked $(d-1)$-sphere can be defined by taking the connected sum of several disjoint copies of the boundary of a $d$-simplex.  Such spheres are never balanced as the graph of a $d$-simplex is the complete graph on $d+1$ vertices.  Instead, a stacked cross-polytopal sphere is defined by taking the connected sum of several disjoint copies of the boundary of a $d$-dimensional cross-polytope.  The boundary of a $d$-dimensional cross-polytope is balanced, and so when taking connected sums, vertices of the same color are identified in order to preserve the balanced condition. We refer to this operation as the \textit{balanced connected sum}.

The inequality \eqref{balancedLBT-inequality} was proved for balanced spheres by Goff, Klee, and Novik \cite[Theorem 5.3]{Goff-Klee-Novik-11} and for balanced manifolds by Browder and Klee \cite[Theorem 4.1]{Browder-Klee-11}.  Klee and Novik extended the inequality  to the family of balanced normal pseudomanifolds  \cite[Theorem 3.4]{Klee-Novik-14} and treated the case of equality \cite[Theorem 4.1]{Klee-Novik-14}.

The following conjecture seems to be a natural extension of Theorem \ref{LBT-Murai}.

\begin{conjecture} {\rm (\cite[Conjecture 3.8]{Klee-Novik-14})\;}
Let $\Delta$ be a balanced normal pseudomanifold of dimension $d-1 \geq 3$.  Then $$2h_2(\Delta) - (d-1)h_1(\Delta) \geq 4 \binom{d }{2}\left(\beta_1(\Delta; \field) - \beta_0(\Delta; \field)\right).$$
\end{conjecture}

The case of equality in this conjecture can be achieved by the complexes in the balanced analog of the Walkup class; that is, the family of complexes defined by beginning with disjoint copies of the boundary of the $d$-dimensional cross-polytope, then applying the operations of balanced connected sums and balanced handle additions.  At the moment, we only know that if $\beta_1(\Delta;\mathbb{Q}) \neq 0$, then $2h_2(\Delta) - (d-1) h_1(\Delta) + 4\binom{d}{2} \beta_0(\Delta) \geq 4\binom{d }{2}$ (see \cite[Theorem 3.9]{Klee-Novik-14}).

The inequality \eqref{balancedLBT-inequality} is equivalent to the inequality $\binom{d}{1} h_2(\Delta) \geq \binom{d}{2} h_1(\Delta)$. This leads to the following balanced version of the Generalized Lower Bound Conjecture.

\begin{conjecture} {\rm (Balanced GLBC \cite[Conjecture 5.6]{Klee-Novik-14})\;}
Let $\Delta$ be a balanced $d$-polytope (or more generally, a balanced $\field$-homology $(d-1)$-sphere).  Then
\begin{equation} \label{balancedGLBC-inequality}
1 = \frac{h_0(\Delta)}{\binom{d}{0}} \leq \frac{h_1(\Delta)}{\binom{d}{1}} \leq \frac{h_2(\Delta)}{\binom{d}{2}} \leq \cdots \leq \frac{h_{\lfloor d/2 \rfloor}(\Delta)}{\binom{d}{\lfloor d/2 \rfloor}}.
\end{equation}
\end{conjecture}

Very recently, Juhnke-Kubitzke and Murai \cite{Kubitzke-Murai-15} established the inequalities \eqref{balancedGLBC-inequality} for all balanced polytopes.  Klee and Novik also proposed a treatment of the case of equality at any step of \eqref{balancedGLBC-inequality}, but it requires additional terminology that we do not wish to introduce here.  See \cite[Section 5.2]{Klee-Novik-14} for further details.

\subsection{Flag complexes}

\begin{definition}
Let $\Delta$ be a simplicial complex.  A subset $S \subseteq V(\Delta)$ is a \textbf{missing face} in $\Delta$ if $S \notin \Delta$ but $T \in \Delta$ for all proper subsets $T \subset S$.  We say that $\Delta$ is \textbf{flag} if all of its missing faces have cardinality $2$.
\end{definition}

A flag simplicial complex is completely determined by its graph: $F \subseteq V(\Delta)$ spans a face in $\Delta$ if and only if all pairs of vertices in $F$ are connected by edges.  For this reason, flag simplicial complexes are also called clique complexes of graphs (or dually, independence complexes of graphs).  The barycentric subdivision of a simplicial complex is another canonical example of a flag complex.

The $h$-vector of a simplicial sphere is palindromic by the Dehn-Sommerville equations.  Therefore, it is natural to expand the $h$-polynomial of a simplicial sphere in a symmetric polynomial basis as $$\sum_{j=0}^dh_j(\Delta)t^j = \sum_{i=0}^{\lfloor d/2 \rfloor} \gamma_i(\Delta)\cdot t^i(1+t)^{d-2i}.$$ The coefficients $\gamma_0, \gamma_1, \ldots, \gamma_{\lfloor d/2 \rfloor}$ are called the \textbf{$\gamma$-numbers} of $\Delta$. The $\gamma$-numbers are weighted alternating sums of the $h$-numbers; in particular, they are integers.

\subsubsection{Lower bounds for flag complexes}

Gal \cite{Gal-05} conjectured the following strengthening of the Generalized Lower Bound Conjecture to the class of flag homology spheres. 

\begin{conjecture}{\rm (Gal's Conjecture \cite[Conjecture 2.1.7]{Gal-05})\;}
If $\Delta$ is a flag homology sphere, then $\gamma_i(\Delta) \geq 0$ for all $i$. 
\end{conjecture}

In the case that $i = \lfloor d/2 \rfloor$, this conjecture is originally due to Charney and Davis \cite{Charney-Davis-95}.  In the language of $h$-vectors, their conjecture states that if $\Delta$ is a flag homology sphere of dimension $2k-1$, then $$(-1)^k\sum_{j=0}^{2k}(-1)^jh_j(\Delta) \geq 0.$$  In the even-dimensional case, the sum on the left-hand side of the above inequality always equals zero by the Dehn-Sommerville equations.

Presently, Gal's conjecture is known to hold for all flag $\mathbb{Q}$-homology $3$-spheres \cite{Davis-Okun-01} and also for barycentric subdivisions of regular CW spheres \cite{Karu-06}.

If Gal's conjecture is true, then the inequalities of Equation \eqref{balancedGLBC-inequality} (the balanced GLBC inequality) continue to hold for flag homology spheres as follows.  Proving the balanced GLBC is equivalent to verifying that $j \cdot h_j(\Delta) \geq (d-j+1)\cdot h_{j-1}(\Delta)$ for all $j \leq \lfloor d/2 \rfloor$.  And, indeed,
\begin{eqnarray*}
j \cdot h_j(\Delta) &=& \sum_{i=0}^j \gamma_i(\Delta)\cdot j \cdot\binom{d-2i}{j-i} 
\geq \sum_{i=0}^{j-1} \gamma_i(\Delta) \cdot j \cdot \binom{d-2i}{j-i} \qquad (\text{ since } \gamma_j \geq 0) \\
&\geq& \sum_{i=0}^{j-1} \gamma_i(\Delta) \cdot (d-j+1)\binom{d-2i}{j-1-i} 
= (d-j+1) \cdot h_{j-1}(\Delta).
\end{eqnarray*}
The final inequality uses the fact that $j\binom{d-2i}{ j-i} \geq (d-j+1)\binom{d-2i}{j-i-1}$ for $j \leq \lfloor d/2 \rfloor$.

Thus the following conjecture of Nevo (personal communication) may be viewed as a relaxation of Gal's conjecture.

\begin{conjecture}
The cone generated by the set of $h$-vectors of balanced $d$-polytopes (or balanced homology $(d-1)$-spheres) contains the cone generated by the set of $h$-vectors of flag $d$-polytopes (or flag homology $(d-1)$-spheres).  
\end{conjecture}
At present, this conjecture is completely open, even in the polytopal case.

\subsubsection{Attempting to characterize $f$-vectors of flag complexes}

While Gal's conjecture can be viewed as a flag analog of the GLBC, the following increasingly strong conjectures of Nevo and Petersen \cite{Nevo-Petersen-11} can be viewed as flag analogs of the necessity part of the $g$-conjecture. 

\begin{conjecture} \label{conj:Nevo-Petersen}
Let $\Delta$ be a flag homology sphere.  Then there exists a simplicial complex $\Gamma$ such that
\begin{enumerate}
\item \cite[Conjecture 1.4]{Nevo-Petersen-11} $\gamma(\Delta) = f(\Gamma)$; 
\item  \cite[Conjecture 6.3]{Nevo-Petersen-11} $\Gamma$ is balanced and $\gamma(\Delta) = f(\Gamma)$; 
\item  \cite[Problem 6.4]{Nevo-Petersen-11} $\Gamma$ is flag and $\gamma(\Delta) = f(\Gamma)$.
\end{enumerate}
\end{conjecture}

The reason that the third part of this conjecture is stronger than the second is due to a result of Frohmader \cite{Frohmader-08}, which states that for every flag simplicial complex, there exists a balanced simplicial complex with the same $f$-vector.

Nevo and Petersen showed that the version of this conjecture asking that $\Gamma$ should be balanced holds through dimension four \cite{Nevo-Petersen-11}.  Nevo, Petersen, and Tenner \cite{Nevo-Petersen-Tenner-11} proved that the same version holds when $\Delta$ is the barycentric subdivision of a homology sphere. Furthermore, Nevo and Petersen verified that the version of this conjecture asking that $\Gamma$  should be flag holds for Coxeter complexes and for flag $(d-1)$-spheres with at most $2d+3$ vertices (i.e., when $\gamma_1 \leq 3$).  

\subsubsection{Upper bounds on flag complexes}

If the second part of Conjecture \ref{conj:Nevo-Petersen} is true, then standard facts about Tur\'an graphs would imply upper bounds on the $\gamma$-numbers of flag spheres. The following conjecture was proposed by Lutz and Nevo \cite[Conjecture 6.3]{Lutz-Nevo-14}.

\begin{conjecture}
Let $d=2k$ be even and suppose $d \geq 4$.  Let $\Delta$ be a flag homology $(d-1)$-sphere on $n$ vertices.  Then $$\gamma_i(\Delta) \leq \gamma_i(\mathcal{J}_{k}(n)),$$ where $\mathcal{J}_k(n)$ is the flag sphere on $n$ vertices obtained as the join of $k$ (graph) cycles, each of which has $\lfloor \frac{n}{k} \rfloor$ or $\lceil \frac{n}{k} \rceil$ vertices.  Moreover, equality holds for some $2 \leq i \leq k$ if and only if $\Delta = \mathcal{J}_k(n)$.
\end{conjecture}

Very recently, Adamaszek and Hladk\'y announced a proof of this conjecture for complexes with a sufficient number of vertices.  For $n\gg 0$, they proved the following stronger result: 

\begin{theorem} {\rm \cite{Adamaszek-Hladky-15}\;}
For every even $d=2k$ with $d \geq 4$, there is a constant $N_0$ for which the following holds: if $M$ is a flag $\mathbb{Z}$-homology manifold of dimension $d-1$ with $n \geq N_0$ vertices, then
\begin{itemize}
\item $f_i(M) \leq f_i(\mathcal{J}_{k}(n))$ for all $1 \leq i \leq d-1$; 
\item $h_i(M) \leq h_i(\mathcal{J}_{k}(n))$ for all $2 \leq i \leq d-2$; 
\item $g_i(M) \leq g_i(\mathcal{J}_{k}(n))$ for all $2 \leq i \leq k$; and 
\item $\gamma_i(M) \leq \gamma_i(\mathcal{J}_{k}(n))$ for all $2 \leq i \leq k$.
\end{itemize}
Moreover, if equality holds in any of these cases, then $M$ is isomorphic to $\mathcal{J}_{k}$.
\end{theorem}

\section{Tools and techniques}   \label{sect:tools}

In this section we will introduce some of the key tools used in the proofs of the results discussed in Sections \ref{Classical Results}-\ref{balanced-and-flag}.

\subsection{The Stanley-Reisner ring}  \label{subsection:SR-ring}

Throughout this section, we assume that $\Delta$ is a $(d-1)$-dimensional simplicial complex with $n$ vertices and  $\field$ is an infinite field of an arbitrary characteristic.  We identify the vertices of $\Delta$ with $[n] = \{1,2,\ldots,n\}$. The standard reference to this material is Stanley's book \cite{Stanley-96}.

We consider a polynomial ring $S:=\field[x_1,x_2,\ldots,x_n]$ with one generator for each vertex in $\Delta$.  The \textbf{Stanley-Reisner ideal} of $\Delta$ is $$I_{\Delta}:= \left(x_{i_1}x_{i_2}\cdots x_{i_k} \ : \ \{i_1,i_2,\ldots,i_k\} \notin \Delta \right).$$
The \textbf{Stanley-Reisner ring} (or face ring) of $\Delta$ is the quotient $\field[\Delta]:= S/I_{\Delta}$.  Since $I_{\Delta}$ is a monomial ideal, the quotient ring $\field[\Delta]$ is graded by degree.  The definition of $I_{\Delta}$ ensures that, as a $\field$-vector space, each graded piece of $\field[\Delta]$, denoted $\field[\Delta]_i$, has a natural basis of monomials whose supports correspond to faces of $\Delta$.  The following result shows that the $h$-numbers arise naturally in this algebraic setting. 

\begin{theorem} {\rm \cite[Proposition 3.2]{Stanley-75} \;}
The Poincar\'e-Hilbert series of $\field[\Delta]$ can be expressed as a rational function of the form 
\begin{equation} \label{Hilb-series-face-ring}
F(\field[\Delta];t):= \sum_{i \geq 0} \dim_{\field}(\field[\Delta])_i \cdot t^i = \frac{h_0(\Delta) + h_1(\Delta) \cdot t + \cdots + h_d(\Delta)\cdot t^d}{(1-t)^d}.
\end{equation}
As a consequence, the Krull dimension of $\field[\Delta]$ is $d = \dim(\Delta)+1$.
\end{theorem}

\begin{definition}
A sequence of linear forms, $\theta_1, \theta_2, \ldots, \theta_d \in S$ is called a \textbf{linear system of parameters} (or l.s.o.p.) if the Artinian reduction $\field[\Delta]/(\Theta)$ is a finite-dimensional $\field$-vector space; here $(\Theta):=(\theta_1,\ldots,\theta_d)$.
\end{definition}

It is a general fact that if $\field$ is an infinite field, then $\field[\Delta]$ admits an l.s.o.p.; moreover, any choice of generic linear forms $\theta_1,\ldots,\theta_d \in S$ is an l.s.o.p.  If $\theta_1, \ldots, \theta_d$ is an l.s.o.p. for $\Delta$, we write $\field(\Delta; \Theta)=\field(\Delta):= \field[\Delta]/(\Theta)$.  

Stanley \cite{Stanley-75} used a result of Reisner \cite{Reisner-76} to build a first bridge between the algebraic structure of the Stanley-Reisner ring and the combinatorial structure of sphere triangulations.
\begin{theorem} \label{CM-h-vec-is-M-vec}  {\rm \cite[Section 4]{Stanley-75}\:}
Let $\Delta$ be a $\field$-homology sphere or $\field$-homology ball, and let $\theta_1, \ldots, \theta_d$ be any l.s.o.p. for $\field[\Delta]$.  Then the Poincar\'e-Hilbert series of $\field(\Delta)$ is
\begin{equation} \label{Hilb-series-artinian-reduction}
F(\field(\Delta);t) = \sum_{j=0}^d h_j(\Delta)\cdot t^j.
\end{equation}
In other words, $\dim_{\field}\field(\Delta)_j = h_j(\Delta)$ for all $j$.
\end{theorem}
The key connection between Equations \eqref{Hilb-series-face-ring} and \eqref{Hilb-series-artinian-reduction} is that, for homology spheres and homology balls, the face ring $\field[\Delta]$ is free as a $\field[\theta_1, \ldots, \theta_d]$-module. (Simplicial complexes with this property are called \textbf{Cohen-Macaulay} complexes.) Equation \eqref{Hilb-series-artinian-reduction} then follows from the fact that the Poincar\'e-Hilbert series of $\field[\theta_1, \ldots, \theta_d]$ is $\frac{1}{(1-t)^d}$.

One immediate consequence of Theorem \ref{CM-h-vec-is-M-vec} is  that the $h$-numbers of a simplicial homology sphere or ball are nonnegative.  In fact, a classical theorem of Macaulay \cite{Macaulay-26} can be used to get much stronger inequalities from Theorem \ref{CM-h-vec-is-M-vec}. First, we require two definitions.  

A \textbf{multicomplex}, $\mathcal{M}$, is a family of monomials that is closed under divisibility; i.e., if $\mu \in \mathcal{M}$ and $\nu$ divides $\mu$, then $\nu \in \mathcal{M}$.  We write $\mathcal{M}_i := \{\mu \in \mathcal{M}\ : \ \deg(\mu) = i\}$.  

Next, it is easy to verify that for all natural numbers $m$ and $i$, there exists a unique decomposition of the form $$m = \binom{a_i}{i} + \binom{a_{i-1}}{i-1} + \cdots + \binom{a_j}{j}, \quad \mbox{ where $a_i > a_{i-1} > \cdots > a_j \geq j > 0$}.$$  Define $$m^{\langle i \rangle}:= \binom{a_i + 1}{ i+1} + \binom{a_{i-1}+1}{i} + \cdots + \binom{a_j + 1}{j+1}.$$  A sequence of integers $(1, F_1, F_2, \ldots)$ is an \textbf{$M$-sequence} if $0 \leq F_{i+1} \leq F_i ^{\langle i \rangle}$ for all $i$.

\begin{theorem} {\rm (Macaulay \cite{Macaulay-26})\;} \label{macaulays-thm}
The following conditions are equivalent.
\begin{enumerate}
\item The sequence $(1,F_1, F_2, \ldots)$ of nonnegative integers is an $M$-sequence.
\item There exists a multicomplex $\mathcal{M}$ such that $|\mathcal{M}_i| = F_i$ for all $i$.   
\item There exists a standard graded $\field$-algebra $R$ such that $\dim_{\field}R_i =  F_i$ for all $i$.
\end{enumerate}
\end{theorem}

\begin{corollary} {\rm \cite[Corollary 5.2]{Stanley-75}\;} \label{M-seq-CM}
If $\Delta$ is a homology sphere or homology ball, then $h(\Delta)$ is an $M$-sequence.
\end{corollary}

\subsection{Lefschetz properties and the $g$-theorem}\label{sec:tools-lefschetz}

When $\Delta$ is a homology sphere, Theorem \ref{CM-h-vec-is-M-vec} and the Dehn-Sommerville equations imply that $$\dim_{\field}\field(\Delta)_j = h_j(\Delta) = h_{d-j}(\Delta) = \dim_{\field}\field(\Delta)_{d-j}.$$  The inner equality is a statement that is both combinatorial and topological in nature, while the outer equalities are algebraic.  As such, it is natural to ask whether the equality $\dim_{\field}\field(\Delta)_j = \dim_{\field}\field(\Delta)_{d-j}$ can be seen as a sort of algebraic duality.  

\begin{definition} \quad 
\begin{enumerate}
\item Let $\Delta$ be a homology $(d-1)$-sphere.  We say that $\field[\Delta]$ (or simply $\Delta$) has the \textbf{strong Lefschetz property} (SLP) if there exists an l.s.o.p. $\theta_1, \ldots, \theta_d$ for $\field[\Delta]$ and a linear form $\omega$ such that the multiplication map $$\times \omega^{d-2j}: \field(\Delta;\Theta)_j \rightarrow \field(\Delta;\Theta)_{d-j}$$ is an isomorphism for all $j \leq \lfloor d/2 \rfloor$.
\item Let $\Delta$ be a homology $(d-1)$-sphere or ball. We say that $\field[\Delta]$ has the \textbf{weak Lefschetz property} (WLP) if there exists an l.s.o.p. $\theta_1, \ldots, \theta_d$ for $\field[\Delta]$ and a linear form $\omega$ such that the multiplication map $\times \omega: \field(\Delta; \Theta)_i \rightarrow \field(\Delta; \Theta)_{i+1}$ is either injective or surjective for all $i<d$.
\end{enumerate}
\end{definition}

A few comments are in order. First, it is well-known that if $\field[\Delta]$ has the SLP, then the desired multiplication map is an isomorphism for any generic choice of $\theta_1, \ldots, \theta_d$ and $\omega$ (see \cite[Proposition 3.6]{Swartz-06}). Second, it follows easily from the definition that if $\field[\Delta]$ has the SLP, then it also has the WLP. Third, if $\Delta$ is a homology sphere or ball with the WLP, then $h(\Delta)$ is unimodal.  Indeed, suppose the multiplication map $\times \omega: \field(\Delta; \Theta)_s \rightarrow \field(\Delta;\Theta)_{s+1}$ is surjective for some $s$.  Then $\left(\field(\Delta;\Theta)/(\omega)\right)_{s+1}=0$, so $\left(\field(\Delta;\Theta)/(\omega)\right)_{i}=0$ for all $i\geq s+1$, and hence $\times \omega: \field(\Delta; \Theta)_{s+1} \rightarrow \field(\Delta; \Theta)_{s+2}$ is also surjective. In particular, if $\Delta$ is a homology sphere with the WLP, then the Dehn-Sommerville equations imply that $h(\Delta)$ is unimodal with its peak at $h_{\lfloor d/2 \rfloor}(\Delta)$.

In fact, the necessity of the conditions of the $g$-theorem for simplicial polytopes is an easy consequence of the (strong or weak) Lefschetz property: the injectivity of the map $\times \omega: \field(P)_i \rightarrow \field(P)_{i+1}$ for all $i \leq \lfloor d/2 \rfloor $ implies that the quotient ring $\field(P)/(\omega)$ is a standard graded $\field$-algebra whose $i$-th graded piece has dimension $g_i(P) = h_i(P) - h_{i-1}(P)$ for all $i\leq \lfloor d/2\rfloor$.  Thus the $g$-vector of $P$ is an $M$-sequence thanks to the following theorem. 

\begin{theorem} {\rm (Stanley \cite{Stanley-80}, McMullen \cite{McMullen-93, McMullen-96})\;} \label{strongLefschetz} A simplicial $d$-polytope has the strong Lefschetz property over $\mathbb{Q}$.
\end{theorem}

More generally, to prove the $g$-conjecture for homology spheres, it would be sufficient to verify that any homology sphere $\Delta$ has the WLP. In fact, it is even sufficient to prove that there exists a one-form $\omega$ such that $\times\omega \ : \ \field(\Delta)_{\lfloor d/2 \rfloor} \rightarrow \field(\Delta)_{\lfloor d/2 \rfloor + 1}$ is surjective (see  \cite[Theorem 3.2]{Novik-Swartz-09:Gorenstein}).

\subsection{From spheres to manifolds: think global, act local}

The main philosophy of many papers on the $f$-numbers of manifolds and pseudomanifolds, is that even though such a space can be very complicated globally, its local structure controls the $f$-numbers (and other invariants) to a large degree. In this section, we demonstrate several examples of this principle.

One such example is given by the $\mu$-numbers of Section \ref{section:sigma-and-mu}. These numbers are defined as sums of certain invariants of the vertex links, yet they provide some information on the Betti numbers of the original complex (see Theorem \ref{mu-vs-beta}).

In a similar spirit, if $\Delta$ is a pure $(d-1)$-dimensional simplicial complex, then the \textbf{short simplicial $h$-vector} of $\Delta$,  $\hat{h}(\Delta)=(\hat{h}_0(\Delta), \hat{h}_1(\Delta),\ldots, \hat{h}_{d-1}(\Delta))$, is defined by $\hat{h}_i(\Delta):=\sum_{v\in V} h_i(\lk_\Delta(v))$ for $i=0,1,\ldots, d-1$. This vector was introduced in \cite{Hersh-Novik-02}. It has a few convenient and easy-to-check properties:

\begin{proposition}  \label{short-h}
The $f$-numbers of a simplicial complex $\Delta$ are non-negative linear combinations of $\hat{h}$-numbers of $\Delta$. Moreover, if $\Delta$ is a homology $(d-1)$-manifold, then $\hat{h}_i(\Delta)=\hat{h}_{d-1-i}(\Delta)$ for all $i=0,1,\ldots,d-1$.
\end{proposition}

One immediate consequence of these two simple facts, together with Stanley's Upper Bound Theorem for spheres (the $h$-vector version), is the proof of the $f$-vector version of the Upper Bound Theorem for all odd-dimensional manifolds (see \cite{Hersh-Novik-02} for a stronger result).

Yet another example of this phenomenon is the following deep observation of Swartz \cite[Theorem 4.26]{Swartz-09} on the existence of Lefschetz-like elements. We denote by $L^i_s(\Delta)$ the set of all pairs $(\omega,\Theta)$ such that $\omega$ is a linear form, $\Theta$ is an l.s.o.p., and the multiplication map $\times \omega: \field(\Delta;\Theta)_i\to \field(\Delta;\Theta)_{i+1}$ is a surjection.

\begin{theorem}  \label{L^i_s}
Let $\Delta$ be a $\field$-homology manifold with $\field$ an infinite field. If  $L^i_s(\lk_{\Delta}(v))$ is nonempty for (almost) all of the vertices $v$ of $\Delta$, then $L^{i+1}_s(\Delta)$ is nonempty.
\end{theorem}

\subsection{Schenzel's formula and socles}  \label{Schenzel-formula}

Many results on the face numbers of manifolds bound the $f$-numbers or $h$-numbers of a simplicial manifold in terms of its Betti numbers.  The proofs of most of these results rely on the following theorem of Schenzel \cite{Schenzel-81}, which relates the face ring of a simplicial manifold to its topological Betti numbers.

\begin{theorem} {\rm \cite{Schenzel-81}\;} \label{thm:Schenzel1}
Let $\Delta$ be a $\field$-homology manifold (with or without boundary) of dimension $d-1$, and let $\theta_1, \ldots, \theta_d$ be an l.s.o.p. for $\field[\Delta]$.  Then
\begin{equation} \label{Schenzels-formula}
\dim_{\field} \field(\Delta)_j = h_j(\Delta) + \binom{d}{j} \sum_{i=1}^{j-1} (-1)^{j-i-1}\beta_{i-1}(\Delta; \field).
\end{equation}
\end{theorem}

In light of this result, we define the $h'$-numbers of $\Delta$ as $$h'_j(\Delta) := h_j(\Delta) + \binom{d}{ j} \sum_{i=1}^{j-1} (-1)^{j-i-1}\beta_{i-1}(\Delta; \field).$$  Note that if $\Delta$ is a $\field$-homology sphere or ball, then $h'(\Delta)=h(\Delta)$. In particular, Theorem \ref{thm:Schenzel1} provides a vast generalization of Theorem \ref{CM-h-vec-is-M-vec}. Also since $\field(\Delta)$ is a standard graded algebra, Schenzel's formula, together with Macaulay's Theorem (Theorem \ref{macaulays-thm}) imply that if $\Delta$ is a $\field$-homology manifold, then $h'(\Delta):=(h_0'(\Delta),\ldots, h'_d(\Delta))$ is an $M$-sequence (cf.~Corollary \ref{M-seq-CM}). 
In order to derive Equation \eqref{Schenzels-formula}, Schenzel \cite{Schenzel-81} proved the following result. 

\begin{theorem} {\rm \cite{Schenzel-81}\;} \label{thm:Schenzel2}
Let $\Delta$ be a $\field$-homology manifold (with or without boundary) of dimension $d-1$, and let $\theta_1, \ldots, \theta_d$ be an l.s.o.p. for $\field[\Delta]$.  For all $1 \leq s \leq d$ and $0 \leq j \leq d$, the kernel of the multiplication map $$\times \theta_s : \left(\field[\Delta]/(\theta_1,\ldots,\theta_{s-1})\right)_{j} \rightarrow \left(\field[\Delta]/(\theta_1,\ldots,\theta_{s-1})\right)_{j+1}$$ has dimension $\binom{s-1}{j} \beta_{j-1}(\Delta; \field)$.  
\end{theorem}

Recall that if $R$ is a commutative ring and $I, J$ are ideals of $R$, then the \textit{ideal quotient} $I:J$ is defined as $I:J = \{x \in R\ : \ xJ \subseteq I\}$.  
\begin{definition}
The \textbf{socle} of $\field(\Delta)$ is the ideal $$ \Soc(\field(\Delta)) = 0: \mathfrak{m} = \{f \in \field(\Delta)\ : \ x_i \cdot f = 0 \text{ for all } i\}.$$  Here we use $\mathfrak{m}$ to denote the maximal ideal $(x_1,\ldots,x_n)$.
\end{definition}

Assume $\Delta$ is a $\field$-homology manifold and $\omega$ is a generic linear form. In view of Theorem \ref{thm:Schenzel2}, one natural question to ask is what can be said about the dimension of the kernel of $\times \omega: 
\field(\Delta)_j\to\field(\Delta)_{j+1}$ or at least about the dimension of the socle of $\field(\Delta)$. To this end, we have:

\begin{theorem} {\rm \cite[Theorem 2.2]{Novik-Swartz-09:Socles}\;} \label{thm:socles}
Let $\Delta$ be a $\field$-homology manifold, with or without boundary.  Then $$\dim_{\field}\Soc(\field(\Delta))_j \geq \binom{d}{ j} \beta_{j-1}(\Delta; \field).$$
\end{theorem}

This theorem is due to Novik and Swartz; its proof relies upon the local cohomology module of $\field[\Delta]$ and the following result of Hochster, which describes the structure of $H^i_{\mathfrak{m}}(\field[\Delta])$ as a $\field$-vector space in terms of the simplicial (co)homology of $\Delta$ and its links.

\begin{theorem}  \label{Hochster's-thm}
Let $\Delta$ be a simplicial complex on $n$ vertices.  Let $\mathbf{a} = ( a_1, \ldots, a_n) \in \mathbb{Z}^{n}$ and let $F = \{j \in [n] \ : \ a_j \neq 0\}$.  Then for any $i \geq 0$, the $i^{\text{th}}$ local cohomology module of $\field[\Delta]$ satisfies
\begin{equation*}
\dim_{\field}H^i_{\mathfrak{m}}(\field[\Delta])_{\mathbf{a}} = 
\begin{cases}
0 & \text{if } F \notin \Delta \text{ or } a_j > 0 \text{ for some } j,  \\
\beta_{i-|F|-1}(\lk_{\Delta}(F); \field) & \text{otherwise}.
\end{cases}
\end{equation*}
\end{theorem}

Let $\Delta$ be a $\field$-homology manifold (with or without boundary). As we noted earlier, the vector $h'(\Delta)$ is then an $M$-sequence. One easy but important consequence of Theorem \ref{thm:socles} is the  following strengthening of this fact: 
\begin{proposition} \label{str-hprime-vector} If $\Delta$ is a $\field$-homology manifold (with or without boundary), then for all $j\leq d-1$, the numbers 
\[
\left(h'_0(\Delta), h'_1(\Delta), \ldots, h'_{j-1}(\Delta), h'_j(\Delta) - \binom{d}{j}\beta_{j-1}(\Delta; \field),h'_{j+1}(\Delta)\right)\]
form an $M$-sequence.
\end{proposition}

This proposition, in turn, leads to the proof of several theorems discussed in Section \ref{sect:manifolds}, among them (i) the Upper Bound Theorem for various classes of homology manifolds \cite{Novik-98} (for this result even a certain weaker version of Proposition \ref{str-hprime-vector} is enough), and (ii) K\"uhnel's conjecture on the Euler characteristic of even-dimensional $\field$-homology manifolds (see \cite[Theorem 4.4]{Novik-Swartz-09:Socles} for further details). Theorem \ref{thm:socles} along with Theorem \ref{L^i_s} also provides a key ingredient of the proof of part (1) of Theorem \ref{LBT-Murai} for orientable $\field$-homology manifolds (see \cite[Theorem 5.2]{Novik-Swartz-09:Socles}).

It is worth mentioning that Theorem \ref{thm:socles} can be stated more generally for Buchsbaum graded modules of Krull dimension $d$ over $\field[x_1, \ldots, x_n]$ (\cite[Theorem 2.2]{Novik-Swartz-09:Socles}).  For Buchsbaum modules over local Noetherian rings, a similar result was established by Goto \cite[Proposition 3.6]{Goto-83}.

\subsection{The Gorenstein property}
A graded $\field$-algebra of Krull dimension zero is called \textit{Gorenstein} if its socle is a $1$-dimensional $\field$-vector space (see \cite[p.~50]{Stanley-96} for many other equivalent definitions). If $\Delta$ is a $\field$-homology sphere, then $\field(\Delta)$ is a Gorenstein ring \cite[Theorem II.5.1]{Stanley-96}. One major consequence of being Gorenstein is the symmetry of the Hilbert function. This naturally leads to the question of how far is $\field(\Delta)$ from being  Gorenstein  if $\Delta$ is a homology manifold other than a sphere. The answer due to Novik and Swartz turns out to be surprisingly simple. We require one definition and a bit of discussion.

\begin{definition} Let $\Delta$ be a $\field$-homology $(d-1)$-manifold (with or without boundary). Define 
\[
\overline{\field(\Delta)}:=\field(\Delta)/\bigoplus_{i=1}^{d-1}\Soc(\field(\Delta))_i.
\]
\end{definition}

In the above formula, $\Soc(\field(\Delta))_d$ is not included in the direct sum because that graded piece of the socle is all of $\field(\Delta)_d$.  Further, by the definition of the socle, each individual graded piece $\Soc(\field(\Delta))_i$ is an ideal of $\field(\Delta)$ on its own, and so $\overline{\field(\Delta)}$ is a quotient of $\field(\Delta)$ by an ideal.  

Recall also the definition of the $h''$-numbers (see Definition \ref{h-doubleprime-nums}). Observe that Schenzel's formula (Theorem \ref{thm:Schenzel1}), and the socle theorem (Theorem \ref{thm:socles}) imply that for a $\field$-homology $(d-1)$-manifold $\Delta$, $\dim_{\field} \overline{\field(\Delta)}_i \leq h''_i(\Delta)$. 

\begin{theorem}  \label{manif->Gor}
Let $\Delta$ be a $(d-1)$-dimensional, orientable, $\field$-homology manifold without boundary. Then
\begin{enumerate}
\item {\rm \cite[Theorem 1.3]{Novik-Swartz-09:Gorenstein}\;} 
$\dim_{\field} \Soc(\field(\Delta))_i=\binom{d}{i}\beta_{i-1}(\Delta)$ for all $i\leq d$. In particular, $\dim_{\field} \overline{\field(\Delta)}_i = h''_i(\Delta)$.
\item {\rm \cite[Theorem 1.4]{Novik-Swartz-09:Gorenstein}\;} 
Moreover, if $\Delta$ is connected, then $\overline{\field(\Delta)}$ is Gorenstein.
\end{enumerate}
\end{theorem}

The major ingredient in the proof of Theorem \ref{manif->Gor} is Gr\"abe's generalization \cite{Grabe-84} of Hochster's theorem (Theorem \ref{Hochster's-thm}) that provides an explicit description of the structure of the local cohomology module $H^i_{\mathfrak{m}}(\field[\Delta])$ as a $\field[\Delta]$-module in terms of the simplicial (co)homology of the links of faces of $\Delta$ and maps between them. 

In turn, Theorem \ref{manif->Gor} is one of the main ingredients in the proof of Theorem \ref{h-double-prime}. Indeed, along with some standard properties of Gorenstein rings (see for instance \cite[Theorems I.12.5 and I.12.10]{Stanley-96}), Theorem \ref{manif->Gor} implies that if $\Delta$ is a $(d-1)$-dimensional connected, orientable, $\field$-homology manifold, then 
\begin{equation}  \label{nat-isom}
\Hom_\field (\overline{\field(\Delta)}_i, \field) \mbox{ is naturally isomorphic to } \overline{\field(\Delta)}_{d-i}.
\end{equation}
 In particular, it follows that $h''_i(\Delta)=h''_{d-i}(\Delta)$. A much more substantial consequence of \eqref{nat-isom} along with Theorem \ref{L^i_s}  is the injectivity of the map $\times \omega: \overline{\field(\Delta)}_i \to \overline{\field(\Delta)}_{i+1}$ for any $i<\lfloor d/2\rfloor$ and a generic linear form $\omega$, under an additional assumption that (most of) vertices of $\Delta$ have the WLP (see proof of \cite[Theorem 3.2]{Novik-Swartz-09:Gorenstein}). This injectivity, in turn, yields part (3) of Theorem \ref{h-double-prime}. In fact, Theorems \ref{L^i_s} and \ref{manif->Gor} are also among the key ingredients in the proof of Theorem \ref{tilde-g}, which gives a strengthening of Theorem \ref{h-double-prime}.

\subsection{Generic initial ideals}

Generic initial ideals can be traced at least to Hartshorne \cite{Hartshorne}. They entered the study of the face numbers of simplicial complexes via Kalai's  theory of algebraic shifting \cite{Kalai-84, Bjorner-Kalai-88, Kalai-02}. Since then they have become an indispensable tool, see for instance \cite{Babson-Nevo, Murai-07, Murai-10:Shifting, Novik-98} and many other papers. Perhaps the most spectacular recent application of generic initial ideals is due to Murai and Nevo \cite{Murai-Nevo-13} in their proof of the equality part of the GLBT for simplicial polytopes. 

Below we sketch the definition of generic initial ideals and discuss some of their properties. We order the variables of $S=\field[x_1,\ldots,x_n]$ by $x_1>x_2>\cdots>x_n$, and let $\succ$ be the rev-lex order on the monomials in $x_1,\ldots, x_n$ of the same degree. For instance, in degree 2 we have: $x_1^2 \succ x_1x_2 \succ x_2^2 \succ x_1x_3 \succ x_2x_3 \succ x_3^2 \succ \cdots$. Furthermore, we refine this partial order by declaring that the lower degree monomials are larger than the higher degree monomials. 

For a polynomial $A$ in $S$, denote by $\Init_\succ(A)$ the rev-lex largest monomial that appears in $A$ with a non-zero coefficient.  If $I$ is a homogeneous ideal, then the \textbf{rev-lex initial ideal of $I$}, $\Init_\succ(I)$, is defined as $\Span_\field\{\Init_\succ(A):  A\in I\}$. It is not hard to see that $\Init_\succ(I)$ is in fact an ideal.

Let $\phi=(\phi_{ij})\in\GL_n(\field)$ be a generic matrix.  Then $\phi$ acts on $S_1$ by $\phi(x_i)=\sum_{k=1}^n \phi_{ki}x_k$ and this action can be extended uniquely to a ring automorphism on $S$. If $I$ is a homogeneous ideal, then the \textbf{(rev-lex) generic initial ideal of $I$} is
\[
\Gin_\succ(I):=\Init_\succ(\phi(I)).
\]
It is a result of Galligo (see, for instance, \cite[Theorem 1.27]{Green-98}) that $\Init_\succ(\phi(I))$ is independent of $\phi$ for generic choices of $\phi$, and so $\Gin_\succ(I)$ is well-defined.

One beautiful property of generic initial ideals is that they have a very simple structure. For instance, if $\field$ is a field of characteristic zero, $I$ is a homogeneous ideal in $S$, and $\mu\in \Gin_\succ(I)$ is a monomial, then for any $x_t$ that divides $\mu$ and any $1\leq s<t$, the monomial $\mu \cdot x_s/x_t$ is also an element of $\Gin_\succ(I)$. (This result is also due to Galligo \cite[Theorem 4.2.1]{Herzog-Hibi}.) In fact, the main reason that generic initial ideals are so useful is that while the degeneration of $I$ to $\Gin_\succ(I)$ considerably simplifies the ideal, it preserves many important properties:

\begin{theorem}
Let $I$ be a homogeneous ideal in $S$.
\begin{enumerate}
\item The ring $S/I$ is Cohen-Macaulay if and only if $S/\Gin_\succ(I)$ is.
\item The ring $S/I$ has the WLP if and only if $S/\Gin_\succ(I)$ has the WLP.
\item {\rm (Crystallization Principle)\;} Suppose that $\field$ is a field of characteristic zero and that $I$ is generated by elements of degree $\leq m$. If $\Gin_\succ(I)$ has no minimal generators of degree $m+1$, then $\Gin_\succ(I)$ is generated
by elements of degree $\leq m$.
\end{enumerate}
\end{theorem}

The first part is due to Bayer and Stillman, see \cite[Corollary 4.3.18]{Herzog-Hibi}. (In fact, a remarkable extension of this part, due to Bayer, Charalambous, and S.~Popescu, asserts that even the values of certain graded Betti numbers are preserved under generic initial ideals; see \cite{Bayer-Charalambous-Popescu-99} and \cite{Aramova-Herzog} for more details.) The second part follows from \cite[Lemma 4.3.7]{Herzog-Hibi} and relies on the properties of the rev-lex order. The third part is due to Green \cite[Proposition 2.28]{Green-98}. This part played a crucial role in the Murai-Nevo proof \cite{Murai-Nevo-13} of the equality part of the GLBT.

\subsection{Graded Betti numbers of the Stanley-Reisner ring}  \label{section:gradedBetti}

Let $S = \field[x_1,\ldots,x_n]$ be a polynomial ring.  It is a well-known fact that if $I \subseteq S$ is a homogeneous ideal, then there is a unique graded minimal free resolution of $S/I$ (viewed as an $S$-module) of the following form: 

\begin{equation*} \label{min-free-resolution}
0\rightarrow\bigoplus_{j}S(-j)^{\beta_{\ell,j}} \stackrel{\varphi_{\ell}}{\rightarrow}
\cdots  \stackrel{\varphi_2}{\rightarrow} 
\bigoplus_{j} S(-j)^{\beta_{1,j}} \stackrel{\varphi_1}{\rightarrow}
\bigoplus_{j} S(-j)^{\beta_{0,j}} \stackrel{\varphi_0}{\rightarrow}
S \rightarrow S/I \rightarrow 0,
\end{equation*}
where $\ell \leq n$.  The notation $S(-j)$ denotes a shift in the degree of $S$: \; $\left(S(-j)\right)_k = S_{j+k}$ for all $k$.  The resolution is \textit{graded} because all maps $\varphi_0, \varphi_1, \ldots, \varphi_{\ell}$ preserve the degrees of elements of $S$ (after the grading on $S$ is shifted). 

The numbers $\beta_{i,j}$ are called the \textbf{graded Betti numbers} of $S/I$.  When $I$ is the Stanley-Reisner ideal of a simplicial complex, the following result of Hochster establishes a beautiful relationship between the graded (algebraic) Betti numbers of the minimal free resolution of $S/I$ and the (topological) Betti numbers of the underlying simplicial complex. 

\begin{theorem} {\rm (Hochster's formula)\; }
Let $\Delta$ be a simplicial complex.  Then $$\beta_{i,i+j}(\field[\Delta]) = \sum_{\stackrel{W \subseteq V(\Delta)}{|W| = i+j}} \beta_{j-1}(\Delta_W; \field).$$
\end{theorem}

Kolins \cite[Proposition 8]{Kolins-11} used the graded Betti numbers to place sufficient conditions on whether a homology ball can be split as the boundary connected sum of two other homology balls along a codimension-one face.  This allowed him to show that several integer vectors could not be realized as $h$-vectors of homology balls, despite satisfying very natural algebraic conditions \cite[Example 10]{Kolins-11}.

More recently, Murai \cite{Murai-15} made a tremendous connection between Hochster's formula and the Bagchi-Datta $\sigma$-numbers by noticing that if $\Delta$ has $n$ vertices, then $$\sigma_{j-1}(\Delta;\field) = \sum_{k=j}^n \frac{\beta_{k-j,k}(\field[\Delta])}{\binom{n}{k}}.$$  Murai then used several results from commutative algebra to provide upper bounds on the graded Betti numbers of $\field[\Delta]$ and hence also the $\sigma$-numbers and $\mu$-numbers.  This, along with the results of \cite{Migliore-Nagel-03}, led to Murai's proof of Theorems \ref{Kuhnel-2} and \ref{LBT-Murai}.

\bigskip We are looking forward to many future results on the face numbers of simplicial complexes and new connections between commutative algebra, algebraic topology, and combinatorics! 

\section*{Acknowledgments}

We are grateful to Bhaskar Bagchi, Basudeb Datta, Satoshi Murai, Eran Nevo, and Ed Swartz for many helpful conversations.


%
%
 \bibliographystyle{plain}
 \bibliography{manifolds-biblio}

\begin{thebibliography}{10}

\bibitem{Adamaszek-Hladky-15}
Micha\l{} Adamaszek and Hladk\'y Jan.
\newblock Upper bound theorem for odd-dimensional flag manifolds.
\newblock \texttt{http://arxiv.org/abs/1503.05961}, 2015.

\bibitem{Aramova-Herzog}
Annetta Aramova and J{\"u}rgen Herzog.
\newblock Almost regular sequences and {B}etti numbers.
\newblock {\em Amer. J. Math.}, 122(4):689--719, 2000.

\bibitem{Babson-Nevo}
Eric Babson and Eran Nevo.
\newblock Lefschetz properties and basic constructions on simplicial spheres.
\newblock {\em J. Algebraic Combin.}, 31(1):111--129, 2010.

\bibitem{Babson-Novik}
Eric Babson and Isabella Novik.
\newblock Face numbers and nongeneric initial ideals.
\newblock {\em Electron. J. Combin.}, 11(2):Research Paper 25, 23, 2004/06.

\bibitem{Bagchi:mu-vector}
Bhaskar Bagchi.
\newblock {The mu vector, Morse inequalities and a generalized lower bound
  theorem for locally tame combinatorial manifolds}.
\newblock \texttt{http://arxiv.org/pdf/1405.5675.pdf}, 2014.

\bibitem{Bagchi-Datta-13}
Bhaskar Bagchi and Basudeb Datta.
\newblock On {$k$}-stellated and {$k$}-stacked spheres.
\newblock {\em Discrete Math.}, 313(20):2318--2329, 2013.

\bibitem{Barnette-73}
David Barnette.
\newblock A proof of the lower bound conjecture for convex polytopes.
\newblock {\em Pacific J. Math.}, 46:349--354, 1973.

\bibitem{Bayer-Charalambous-Popescu-99}
Dave Bayer, Hara Charalambous, and Sorin Popescu.
\newblock Extremal {B}etti numbers and applications to monomial ideals.
\newblock {\em J. Algebra}, 221(2):497--512, 1999.

\bibitem{Benedetti}
Bruno Benedetti.
\newblock {Smoothing discrete Morse theory}.
\newblock {\em Ann. Sc. Norm. Super. Pisa Cl. Sci.}, to appear;
  \texttt{http://arxiv.org/pdf/1212.0885v4.pdf}, 2014.

\bibitem{Billera-Lee-81}
Louis~J. Billera and Carl~W. Lee.
\newblock A proof of the sufficiency of {M}c{M}ullen's conditions for
  {$f$}-vectors of simplicial convex polytopes.
\newblock {\em J. Combin. Theory Ser. A}, 31(3):237--255, 1981.

\bibitem{Bjorner-Kalai-88}
Anders Bj{\"o}rner and Gil Kalai.
\newblock An extended {E}uler-{P}oincar\'e theorem.
\newblock {\em Acta Math.}, 161(3-4):279--303, 1988.

\bibitem{Brehm-Kuhnel-87}
U.~Brehm and W.~K{\"u}hnel.
\newblock Combinatorial manifolds with few vertices.
\newblock {\em Topology}, 26(4):465--473, 1987.

\bibitem{Brenti-Welker-08}
Francesco Brenti and Volkmar Welker.
\newblock {$f$}-vectors of barycentric subdivisions.
\newblock {\em Math. Z.}, 259(4):849--865, 2008.

\bibitem{Browder-Klee-11}
Jonathan Browder and Steven Klee.
\newblock Lower bounds for {B}uchsbaum complexes.
\newblock {\em European J. Combin.}, 32(1):146--153, 2011.

\bibitem{Cannon-79}
J.~W. Cannon.
\newblock Shrinking cell-like decompositions of manifolds. {C}odimension three.
\newblock {\em Ann. of Math. (2)}, 110(1):83--112, 1979.

\bibitem{Charney-Davis-95}
Ruth Charney and Michael Davis.
\newblock The {E}uler characteristic of a nonpositively curved, piecewise
  {E}uclidean manifold.
\newblock {\em Pacific J. Math.}, 171(1):117--137, 1995.

\bibitem{Chestnut-Sapir-Swartz-08}
Jacob Chestnut, Jenya Sapir, and Ed~Swartz.
\newblock Enumerative properties of triangulations of spherical bundles over
  {$S^1$}.
\newblock {\em European J. Combin.}, 29(3):662--671, 2008.

\bibitem{Datta-Murai}
Basedub Datta and Satoshi Murai.
\newblock On stacked triangulated manifolds.
\newblock \texttt{http://arxiv.org/pdf/1407.6767.pdf}, 2014.

\bibitem{Davis-Okun-01}
Michael~W. Davis and Boris Okun.
\newblock Vanishing theorems and conjectures for the {$\ell^2$}-homology of
  right-angled {C}oxeter groups.
\newblock {\em Geom. Topol.}, 5:7--74, 2001.

\bibitem{Fogelsanger-88}
Allen~Lee Fogelsanger.
\newblock {\em The generic rigidity of minimal cycles}.
\newblock ProQuest LLC, Ann Arbor, MI, 1988.
\newblock Thesis (Ph.D.)--Cornell University.

\bibitem{Frohmader-08}
Andrew Frohmader.
\newblock Face vectors of flag complexes.
\newblock {\em Israel J. Math.}, 164:153--164, 2008.

\bibitem{Gal-05}
{\'S}wiatos{\l}aw~R. Gal.
\newblock Real root conjecture fails for five- and higher-dimensional spheres.
\newblock {\em Discrete Comput. Geom.}, 34(2):269--284, 2005.

\bibitem{Goff-Klee-Novik-11}
Michael Goff, Steven Klee, and Isabella Novik.
\newblock Balanced complexes and complexes without large missing faces.
\newblock {\em Ark. Mat.}, 49(2):335--350, 2011.

\bibitem{Goto-83}
Shiro Goto.
\newblock On the associated graded rings of parameter ideals in {B}uchsbaum
  rings.
\newblock {\em J. Algebra}, 85(2):490--534, 1983.

\bibitem{Grabe-84}
Hans-Gert Gr{\"a}be.
\newblock The canonical module of a {S}tanley-{R}eisner ring.
\newblock {\em J. Algebra}, 86(1):272--281, 1984.

\bibitem{Grabe-87}
Hans-Gert Gr{\"a}be.
\newblock Generalized {D}ehn-{S}ommerville equations and an upper bound
  theorem.
\newblock {\em Beitr\"age Algebra Geom.}, 25:47--60, 1987.

\bibitem{Green-98}
Mark~L. Green.
\newblock Generic initial ideals.
\newblock In {\em Six lectures on commutative algebra ({B}ellaterra, 1996)},
  volume 166 of {\em Progr. Math.}, pages 119--186. Birkh\"auser, Basel, 1998.

\bibitem{Hartshorne}
Robin Hartshorne.
\newblock Connectedness of the {H}ilbert scheme.
\newblock {\em Inst. Hautes \'Etudes Sci. Publ. Math.}, 29:5--48, 1966.

\bibitem{Hersh-Novik-02}
Patricia Hersh and Isabella Novik.
\newblock A short simplicial {$h$}-vector and the upper bound theorem.
\newblock {\em Discrete Comput. Geom.}, 28(3):283--289, 2002.

\bibitem{Herzog-Hibi}
J{\"u}rgen Herzog and Takayuki Hibi.
\newblock {\em Monomial ideals}, volume 260 of {\em Graduate Texts in
  Mathematics}.
\newblock Springer-Verlag London, Ltd., London, 2011.

\bibitem{Kubitzke-Murai-15}
Martina Juhnke-Kubitzke and Satoshi Murai.
\newblock Balanced generalized lower bound inequality for simplicial polytopes.
\newblock \texttt{http://arxiv.org/abs/1503.06430}, 2015.

\bibitem{Jungerman-Ringel-80}
M.~Jungerman and G.~Ringel.
\newblock Minimal triangulations on orientable surfaces.
\newblock {\em Acta Math.}, 145(1-2):121--154, 1980.

\bibitem{Kalai-84}
Gil Kalai.
\newblock Characterization of {$f$}-vectors of families of convex sets in
  {${\bf R}^d$}. {I}. {N}ecessity of {E}ckhoff's conditions.
\newblock {\em Israel J. Math.}, 48(2-3):175--195, 1984.

\bibitem{Kalai-87}
Gil Kalai.
\newblock Rigidity and the lower bound theorem. {I}.
\newblock {\em Invent. Math.}, 88(1):125--151, 1987.

\bibitem{Kalai-88}
Gil Kalai.
\newblock Many triangulated spheres.
\newblock {\em Discrete Comput. Geom.}, 3(1):1--14, 1988.

\bibitem{Kalai-02}
Gil Kalai.
\newblock Algebraic shifting.
\newblock In {\em Computational commutative algebra and combinatorics ({O}saka,
  1999)}, volume~33 of {\em Adv. Stud. Pure Math.}, pages 121--163. Math. Soc.
  Japan, Tokyo, 2002.

\bibitem{Karu-06}
Kalle Karu.
\newblock The {$cd$}-index of fans and posets.
\newblock {\em Compos. Math.}, 142(3):701--718, 2006.

\bibitem{Klee-09}
Steven Klee.
\newblock The fundamental group of balanced simplicial complexes and posets.
\newblock {\em Electron. J. Combin.}, 16(2, Special volume in honor of Anders
  Bj\"orner):Research Paper 7, 12, 2009.

\bibitem{Klee-Novik-12}
Steven Klee and Isabella Novik.
\newblock Centrally symmetric manifolds with few vertices.
\newblock {\em Adv. Math.}, 229(1):487--500, 2012.

\bibitem{Klee-Novik-14}
Steven Klee and Isabella Novik.
\newblock {L}ower {B}ound {T}heorems and a {G}eneralized {L}ower {B}ound
  {C}onjecture for balanced simplicial complexes.
\newblock \texttt{http://arxiv.org/abs/1409.5094}, 2014.

\bibitem{Klee-64}
Victor Klee.
\newblock A combinatorial analogue of {P}oincar\'e's duality theorem.
\newblock {\em Canad. J. Math.}, 16:517--531, 1964.

\bibitem{Klee-64:vertices}
Victor Klee.
\newblock On the number of vertices of a convex polytope.
\newblock {\em Canad. J. Math.}, 16:701--720, 1964.

\bibitem{Kolins-11}
Samuel Kolins.
\newblock {$f$}-vectors of triangulated balls.
\newblock {\em Discrete Comput. Geom.}, 46(3):427--446, 2011.

\bibitem{Kubitzke-Nevo-09}
Martina Kubitzke and Eran Nevo.
\newblock The {L}efschetz property for barycentric subdivisions of shellable
  complexes.
\newblock {\em Trans. Amer. Math. Soc.}, 361(11):6151--6163, 2009.

\bibitem{Kuhnel-90}
Wolfgang K{\"u}hnel.
\newblock Triangulations of manifolds with few vertices.
\newblock In {\em Advances in differential geometry and topology}, pages
  59--114. World Sci. Publ., Teaneck, NJ, 1990.

\bibitem{Kuhnel-95}
Wolfgang K{\"u}hnel.
\newblock {\em Tight polyhedral submanifolds and tight triangulations}, volume
  1612 of {\em Lecture Notes in Mathematics}.
\newblock Springer-Verlag, Berlin, 1995.

\bibitem{Lee-94}
Carl~W. Lee.
\newblock Generalized stress and motions.
\newblock In {\em Polytopes: abstract, convex and computational ({S}carborough,
  {ON}, 1993)}, volume 440 of {\em NATO Adv. Sci. Inst. Ser. C Math. Phys.
  Sci.}, pages 249--271. Kluwer Acad. Publ., Dordrecht, 1994.

\bibitem{Lee-Schmidt-12}
Carl~W. Lee and Laura Schmidt.
\newblock On the numbers of faces of low-dimensional regular triangulations and
  shellable balls.
\newblock {\em Rocky Mountain J. Math.}, 41(6):1939--1961, 2011.

\bibitem{Lutz-05}
Frank~H Lutz.
\newblock Triangulated manifolds with few vertices: Combinatorial manifolds.
\newblock \texttt{http://arxiv.org/pdf/math/0506372v1.pdf}, 2005.

\bibitem{Lutz-Nevo-14}
Frank~H. Lutz and Eran Nevo.
\newblock Stellar theory for flag complexes.
\newblock \texttt{http://arxiv.org/abs/1302.5197}, 2014.

\bibitem{Lutz-Sulanke-Swartz-09}
Frank~H. Lutz, Thom Sulanke, and Ed~Swartz.
\newblock {$f$}-vectors of 3-manifolds.
\newblock {\em Electron. J. Combin.}, 16(2, Special volume in honor of Anders
  Bj\"orner):Research Paper 13, 33, 2009.

\bibitem{Macaulay-26}
F.~S. Macaulay.
\newblock Some {P}roperties of {E}numeration in the {T}heory of {M}odular
  {S}ystems.
\newblock {\em Proc. London Math. Soc.}, S2-26(1):531, 1926.

\bibitem{Macdonald-71}
I.~G. Macdonald.
\newblock Polynomials associated with finite cell-complexes.
\newblock {\em J. London Math. Soc. (2)}, 4:181--192, 1971.

\bibitem{McMullen-70}
P.~McMullen.
\newblock The maximum numbers of faces of a convex polytope.
\newblock {\em Mathematika}, 17:179--184, 1970.

\bibitem{McMullen-93}
P.~McMullen.
\newblock On simple polytopes.
\newblock {\em Invent. Math.}, 113(2):419--444, 1993.

\bibitem{McMullen-96}
P.~McMullen.
\newblock Weights on polytopes.
\newblock {\em Discrete Comput. Geom.}, 15(4):363--388, 1996.

\bibitem{McMullen-Walkup-71}
P.~McMullen and D.~W. Walkup.
\newblock A generalized lower-bound conjecture for simplicial polytopes.
\newblock {\em Mathematika}, 18:264--273, 1971.

\bibitem{Migliore-Nagel-03}
J.~Migliore and U.~Nagel.
\newblock Reduced arithmetically {G}orenstein schemes and simplicial polytopes
  with maximal {B}etti numbers.
\newblock {\em Adv. Math.}, 180(1):1--63, 2003.

\bibitem{Munkres-84}
James~R. Munkres.
\newblock Topological results in combinatorics.
\newblock {\em Michigan Math. J.}, 31(1):113--128, 1984.

\bibitem{Murai-07}
Satoshi Murai.
\newblock Generic initial ideals and squeezed spheres.
\newblock {\em Adv. Math.}, 214(2):701--729, 2007.

\bibitem{Murai-10:Shifting}
Satoshi Murai.
\newblock Algebraic shifting of strongly edge decomposable spheres.
\newblock {\em J. Combin. Theory Ser. A}, 117(1):1--16, 2010.

\bibitem{Murai-10:Barycentric}
Satoshi Murai.
\newblock On face vectors of barycentric subdivisions of manifolds.
\newblock {\em SIAM J. Discrete Math.}, 24(3):1019--1037, 2010.

\bibitem{Murai-15}
Satoshi Murai.
\newblock Tight combinatorial manifolds and graded betti numbers.
\newblock {\em Collect. Math.}, 2015.
\newblock
  \texttt{http://link.springer.com/article/10.1007\%2Fs13348-015-0137-z}.

\bibitem{Murai-Nevo-13}
Satoshi Murai and Eran Nevo.
\newblock On the generalized lower bound conjecture for polytopes and spheres.
\newblock {\em Acta Math.}, 210(1):185--202, 2013.

\bibitem{Murai-Nevo-14}
Satoshi Murai and Eran Nevo.
\newblock On {$r$}-stacked triangulated manifolds.
\newblock {\em J. Algebraic Combin.}, 39(2):373--388, 2014.

\bibitem{Nevo-Petersen-11}
Eran Nevo and T.~Kyle Petersen.
\newblock On {$\gamma$}-vectors satisfying the {K}ruskal-{K}atona inequalities.
\newblock {\em Discrete Comput. Geom.}, 45(3):503--521, 2011.

\bibitem{Nevo-Petersen-Tenner-11}
Eran Nevo, T.~Kyle Petersen, and Bridget~Eileen Tenner.
\newblock The {$\gamma$}-vector of a barycentric subdivision.
\newblock {\em J. Combin. Theory Ser. A}, 118(4):1364--1380, 2011.

\bibitem{Nevo-Santos-Wilson}
Eran Nevo, Francisco Santos, and Stedman Wilson.
\newblock Many triangulated odd-spheres.
\newblock \texttt{http://arxiv.org/pdf/1408.3501.pdf}, 2014.

\bibitem{Novik-98}
Isabella Novik.
\newblock Upper bound theorems for homology manifolds.
\newblock {\em Israel J. Math.}, 108:45--82, 1998.

\bibitem{Novik-05}
Isabella Novik.
\newblock On face numbers of manifolds with symmetry.
\newblock {\em Adv. Math.}, 192(1):183--208, 2005.

\bibitem{Novik-Swartz-09:DS}
Isabella Novik and Ed~Swartz.
\newblock Applications of {K}lee's {D}ehn-{S}ommerville relations.
\newblock {\em Discrete Comput. Geom.}, 42(2):261--276, 2009.

\bibitem{Novik-Swartz-09:Gorenstein}
Isabella Novik and Ed~Swartz.
\newblock Gorenstein rings through face rings of manifolds.
\newblock {\em Compos. Math.}, 145(4):993--1000, 2009.

\bibitem{Novik-Swartz-09:Socles}
Isabella Novik and Ed~Swartz.
\newblock Socles of {B}uchsbaum modules, complexes and posets.
\newblock {\em Adv. Math.}, 222(6):2059--2084, 2009.

\bibitem{Novik-Swartz-12}
Isabella Novik and Ed~Swartz.
\newblock Face numbers of pseudomanifolds with isolated singularities.
\newblock {\em Math. Scand.}, 110(2):198--222, 2012.

\bibitem{Pfeifle-Ziegler}
Julian Pfeifle and G{\"u}nter~M. Ziegler.
\newblock Many triangulated 3-spheres.
\newblock {\em Math. Ann.}, 330(4):829--837, 2004.

\bibitem{Reisner-76}
Gerald~Allen Reisner.
\newblock Cohen-{M}acaulay quotients of polynomial rings.
\newblock {\em Advances in Math.}, 21(1):30--49, 1976.

\bibitem{Ringel}
Gerhard Ringel.
\newblock Wie man die geschlossenen nichtorientierbaren {F}l\"achen in
  m\"oglichst wenig {D}reiecke zerlegen kann.
\newblock {\em Math. Ann.}, 130:317--326, 1955.

\bibitem{Schenzel-81}
Peter Schenzel.
\newblock On the number of faces of simplicial complexes and the purity of
  {F}robenius.
\newblock {\em Math. Z.}, 178(1):125--142, 1981.

\bibitem{Stanley-75}
Richard~P. Stanley.
\newblock The upper bound conjecture and {C}ohen-{M}acaulay rings.
\newblock {\em Studies in Appl. Math.}, 54(2):135--142, 1975.

\bibitem{Stanley-79}
Richard~P. Stanley.
\newblock Balanced {C}ohen-{M}acaulay complexes.
\newblock {\em Trans. Amer. Math. Soc.}, 249(1):139--157, 1979.

\bibitem{Stanley-80}
Richard~P. Stanley.
\newblock The number of faces of a simplicial convex polytope.
\newblock {\em Adv. in Math.}, 35(3):236--238, 1980.

\bibitem{Stanley-96}
Richard~P. Stanley.
\newblock {\em Combinatorics and commutative algebra}, volume~41 of {\em
  Progress in Mathematics}.
\newblock Birkh\"auser Boston, Inc., Boston, MA, second edition, 1996.

\bibitem{Stanley-14}
Richard~P. Stanley.
\newblock How the upper bound conjecture was proved.
\newblock {\em Ann. Comb.}, 18(3):533--539, 2014.

\bibitem{Swartz-06}
Ed~Swartz.
\newblock {$g$}-elements, finite buildings and higher {C}ohen-{M}acaulay
  connectivity.
\newblock {\em J. Combin. Theory Ser. A}, 113(7):1305--1320, 2006.

\bibitem{Swartz-09}
Ed~Swartz.
\newblock Face enumeration---from spheres to manifolds.
\newblock {\em J. Eur. Math. Soc. (JEMS)}, 11(3):449--485, 2009.

\bibitem{Swartz-14}
Ed~Swartz.
\newblock Thirty-five years and counting.
\newblock \texttt{http://arxiv.org/pdf/1411.0987.pdf}, 2014.

\bibitem{Swartz-Oberwolfach}
Ed~Swartz.
\newblock What's next?
\newblock Oberwolfach Workshop on Geometric and Algebraic Combinatorics,
  Oberwolfach Report available at \texttt{http://www.mfo.de/occasion/1506},
  2015.

\bibitem{Tay-95}
Tiong-Seng Tay.
\newblock Lower-bound theorems for pseudomanifolds.
\newblock {\em Discrete Comput. Geom.}, 13(2):203--216, 1995.

\bibitem{Walkup-70}
David~W. Walkup.
\newblock The lower bound conjecture for {$3$}- and {$4$}-manifolds.
\newblock {\em Acta Math.}, 125:75--107, 1970.

\end{thebibliography}

\end{document}